\documentclass[11pt,reqno,tbtags]{article}
\usepackage{amsmath}
\usepackage{amsthm}
\usepackage{url}

\numberwithin{equation}{section}

\newtheorem{theorem}{Theorem}[section]
\newtheorem*{theorem*}{Theorem}
\newtheorem{corollary}[theorem]{Corollary}
\newtheorem{lemma}[theorem]{Lemma}
\newtheorem{proposition}[theorem]{Proposition}
\newtheorem{claim}[theorem]{Claim}

\theoremstyle{definition}
\newtheorem{definition}[theorem]{Definition}

\theoremstyle{remark}
\newtheorem{remark}[theorem]{Remark}
\newtheorem{example}[theorem]{Example}

\newtheorem{identity}[theorem]{Identity}
\newtheorem*{remark*}{Remark}

\newcommand{\lc}{\left\lceil}

\newcommand{\rc}{\right\rceil}

\newcommand{\EE}{{\bf  E}}

\newcommand{\RR}{{\bf  R}}

\newcommand{\CC}{{\bf  C}}

\newcommand{\jj}{{\bf  j}}
\newcommand{\kk}{{\bf  k}}

\newcommand{\Var}{{\bf Var}}

\newcommand{\mut}{{\tilde{\mu}}}

\newcommand{\Xt}{{\tilde{X}}}

\newcommand{\Ah}{\hat{A}}
\newcommand{\Bh}{\hat{B}}
\newcommand{\bh}{\hat{b}}

\newcommand{\begp}{\begin{proposition}}
\newcommand{\enp}{\end{proposition}}
\newcommand{\begt}{\begin{theorem}}
\newcommand{\ent}{\end{theorem}}
\newcommand{\begl}{\begin{lemma}}
\newcommand{\enl}{\end{lemma}}
\newcommand{\begc}{\begin{corollary}}
\newcommand{\enc}{\end{corollary}}
\newcommand{\begcl}{\begin{claim}}
\newcommand{\encl}{\end{claim}}
\newcommand{\begr}{\begin{remark}}
\newcommand{\enr}{\end{remark}}
\newcommand{\begal}{\begin{algorithm}}
\newcommand{\enal}{\end{algorithm}}
\newcommand{\begd}{\begin{definition}}
\newcommand{\enf}{\end{definition}}
\newcommand{\begx}{\begin{example}}
\newcommand{\enx}{\end{example}}
\newcommand{\bega}{\begin{array}}
\newcommand{\ena}{\end{array}}

\newcommand{\sfrac}[2]{{\textstyle\frac{#1}{#2}}}

\def\rompar(#1){\textup(#1\textup)}    

\newcommand\eps{\varepsilon}

\newcommand\ga{\alpha}
\newcommand\gd{\delta}
\newcommand\gl{\lambda}

\newcommand\gth{\theta}

\newcommand{\refS}[1]{Section~\ref{#1}}
\newcommand{\refT}[1]{Theorem~\ref{#1}}

\newcommand{\refL}[1]{Lemma~\ref{#1}}
\newcommand{\refR}[1]{Remark~\ref{#1}}

\renewcommand\Re{\operatorname{Re}}

\newcommand\nopf{\qed}   

\usepackage[final]{graphicx}
\newcommand{\falling}[2]{#1^{\underline{#2}}}
\newcommand{\rising}[2]{#1^{\overline{#2}}}

\newcommand{\stirlingone}[2]{\genfrac{[}{]}{0pt}{}{#1}{#2}}

\newcommand{\indicator}[1]{\mathbf{1}\!\left({#1}\right)}
\newcommand{\E}[1]{\mathbf{E}\,#1}
\usepackage{amssymb}

\newcommand{\sumstark}{\sideset{}{^*}\sum_{\mathbf{k}}}

\begin{document}

\setcounter{page}{0}
\thispagestyle{empty}

\begin{center}
{\Large \bf Transfer Theorems and Asymptotic Distributional Results for $m$-ary
Search Trees}\\
\normalsize

\vspace{4ex}
{\sc James Allen Fill\footnotemark}\\
\vspace{.1in}
Department of Applied Mathematics and Statistics \\
\vspace{.1in}
The Johns Hopkins University \\
\vspace{.1in}
{\tt jimfill@jhu.edu}\; and\; {\tt http://www.ams.jhu.edu/\~{}fill/} \\
\vspace{.18in}
\textsc{and}\\
\vspace{.18in}
{\sc Nevin Kapur\footnotemark[1]\;\footnotemark[2]}\\
\vspace{.1in}
Department of Computer Science \\
\vspace{.1in}
California Institute of Technology \\
\vspace{.1in}
{\tt nkapur@cs.caltech.edu}\; and\; {\tt
  http://www.cs.caltech.edu/\~{}nkapur/} \\

\end{center}
\vspace{3ex}

\begin{center}
{\sl ABSTRACT} \\
\end{center}

We derive asymptotics of moments and identify limiting distributions,
under the random permutation model on $m$-ary search trees, for
functionals that satisfy recurrence relations of a simple additive
form.  Many important functionals including the space requirement,
internal path length, and the so-called shape functional fall under
this framework. The approach is based on establishing \emph{transfer
  theorems} that link the order of growth of the input into a
particular (deterministic) recurrence to the order of growth of the
output.  The transfer theorems are used in conjunction with the method
of moments to establish limit laws.  It is shown that (i)~for small
toll sequences $(t_n)$ [roughly, $t_n = O(n^{1 / 2})$] we have
asymptotic normality if $m \leq 26$ and typically periodic behavior if
$m \geq 27$; (ii)~for moderate toll sequences [roughly, $t_n =
\omega(n^{1 / 2})$ but $t_n = o(n)$] we have convergence to non-normal
distributions if $m \leq m_0$ (where $m_0 \geq 26$) and typically
periodic behavior if $m \geq m_0 + 1$; and (iii)~for large toll sequences
[roughly, $t_n = \omega(n)$] we have convergence to non-normal
distributions for all values of~$m$.

\bigskip
\bigskip

\begin{small}

\par\noindent
{\em AMS\/} 2000 {\em subject classifications.\/}  Primary 68W40;
secondary 60F05, 68P10, 60E05.
%
\medskip
\par\noindent
{\em Key words and phrases.\/}
Transfer theorems, $m$-ary search trees, additive functionals, random
permutation model, limit distribution, Euler differential equation,
indicial polynomial.

\medskip
\par\noindent
\emph{Date.} Revised January~11, 2004.
\end{small}

\footnotetext[1]{Research for both authors supported by NSF grant
  DMS--9803780 and
  DMS--0104167, and by The Johns Hopkins University's Acheson
  J.~Duncan Fund for the Advancement of Research in Statistics.}

\footnotetext[2]{Research supported by NSF grant 0049092 and carried
  out primarily while this author was affiliated with what is now the
  Department of Applied Mathematics and Statistics at The Johns
  Hopkins University.}

\newpage
\addtolength{\topmargin}{+0.5in}


\section{Background and notation}
\label{sec:background}

We start by giving a brief overview of search trees which are fundamental
data structures in computer science used in searching and sorting.
For integer \( m \geq 2 \), the \( m \)-ary 
search tree, or multiway tree, generalizes the binary search tree.
The quantity \( m \) is called the \emph{branching factor}.
According to~\cite{MR93f:68045}, search trees of branching factors
higher than 2 were first suggested by Muntz and
Uzgalis~\cite{muntz71:_dynam} ``to solve internal memory problems with
large quantities of data.''  For more background we refer the reader
to~ \cite{knuth97,knuth98} and~\cite{MR93f:68045}.

\label{def:marytree} An \emph{\( m \)-ary tree} is a rooted tree
with at most \( m \) ``children'' for each \emph{node (vertex)}, each
child of a node being
distinguished as one of \( m \) possible types.
Recursively expressed, an \( m \)-ary tree either is empty or consists
of a distinguished node (called the \emph{root}) together with an
ordered \( m \)-tuple of \emph{subtrees}, each of which is an \( m
\)-ary tree.

\label{def:marysearchtree} An \emph{\( m \)-ary search tree} is an
\( m \)-ary tree in which each node has the capacity to contain \(
m-1 \) elements of some linearly ordered set, called the set of
\emph{keys}.
In typical implementations of \(m\)-ary search trees, the keys at each
node are stored in increasing order and at each node one has \( m \)
pointers to the subtrees. By spreading the input data in \( m \)
directions instead of only 2, as is the case for a binary search tree,
one seeks to have shorter path lengths and thus quicker searches.

We consider the space of \( m \)-ary search trees on \( n \) keys, and
assume that the keys are linearly ordered.  Hence, without loss
of generality, we can take the set of keys to be \( [n] :=
\{1,2,\ldots,n\} \). We construct an \( m \)-ary search tree from a
sequence $s$ of \( n \) distinct keys in the following way:
\begin{enumerate}
\item If \( n < m \), then all the keys are stored in the root node in
  increasing order.\label{item:4}
\item If \( n \geq m \), then the first \( m-1 \) keys in the sequence
      are\label{item:5}
  stored in the root in increasing order, and the remaining \( n-(m-1)
  \) keys are stored in the subtrees subject to the condition that if
  \( \sigma_1 < \sigma_2 < \cdots < \sigma_{m-1} \) denotes the
  ordered sequence of keys in the root, then the keys in the \( j \)th
  subtree are those that lie between \( \sigma_{j-1} \) and \(
  \sigma_{j} \), where \( \sigma_0 := 0 \) and \( \sigma_{m} := n+1
  \), sequenced as in $s$.
\item All the subtrees are \( m \)-ary search trees that satisfy
  conditions~\ref{item:4},~\ref{item:5}, and~\ref{item:6}.\label{item:6}
\end{enumerate}

It is our goal to study additive functionals (see
Definition~\ref{def:additive-functional}) defined on \( m \)-ary search
trees.  Such functionals represent the cost of divide-and-conquer
algorithms, reflecting the inherent recursive nature of the algorithms.

\begin{figure}[htbp]
\begin{center}
  \centering
  \includegraphics{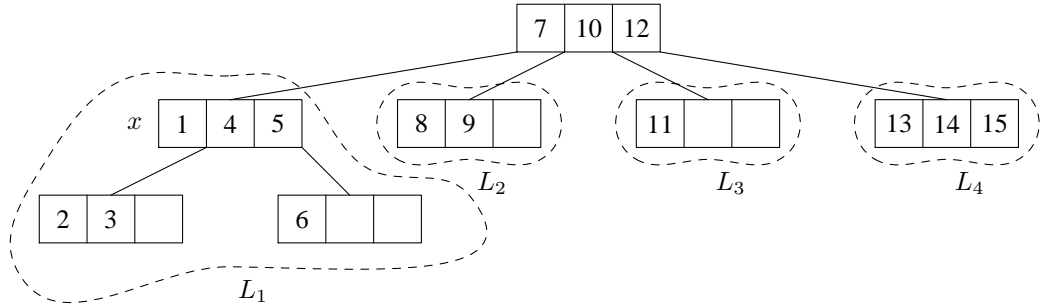}%
  \caption[Example of notation for~$m$-ary search trees.]%
  {Example of notation for a quaternary tree~\(T\):\ here, $|T|=15$
    and $T_x = L_1$.}
  \label{fig:terms}
\end{center}
\end{figure}

Let~\(T\) be an \(m\)-ary search
tree.  We use \(|T|\) to denote the number of keys in \(T\).  Call a
node \emph{full} if it contains \(m-1\) keys.  For \(1 \leq j \leq
m\), let \(L_j(T)\) denote the \(j\)th subtree pendant from the root
of~\(T\).  For a node~$x$ in \( T \),
write \( T_x \) for the subtree of \( T \) consisting of $x$ and its
descendants, with $x$ as root.  This notation is illustrated in
Figure~\ref{fig:terms}.

\begin{definition}
  \label{def:additive-functional} Fix~\( m \geq 2 \).  We will call a
  functional~\( f \) on~\( m \)-ary search trees  an
  \emph{additive tree functional} if it satisfies the recurrence
  \begin{equation}
    \label{eq:2.1}
    f(T) = \sum_{i=1}^m f( L_i(T) ) + t_{|T|},
  \end{equation}
  for any tree~\( T \) with~\( |T| \geq m -1 \). Here~\( (t_n)_{n \geq
  m-1} \) is a given sequence, henceforth called the \emph{toll
  sequence} or \emph{toll function}.
\end{definition}
Note that the recurrence~\eqref{eq:2.1} does not make any reference
to~\( t_n \) for~\( 0 \leq n \leq m-2 \) nor specify~\( f(T) \) for~\(
0 \leq |T| \leq m-2 \).

Several interesting examples can be cast as additive functionals.
\begin{example}
  \label{example:space-requirement} If we specify~\( f(T) \)
  arbitrarily for~\( 0 \leq |T| \leq m-2 \) and take~\( t_n \equiv c
  \) for~\( n \geq m-1 \), we obtain the ``additive functional''
  framework of~\cite[\S3.1]{MR93f:68045}.  (Our definition of an
  additive functional generalizes this notion.)  In particular if we
  define~\( f(\emptyset) := 0 \) and~\( f(T) := 1 \) for the unique \(
  m \)-ary search tree~\( T \) on~\( n \) keys for \( 1 \leq n \leq
  m-2 \) and let \( t_n \equiv 1 \) for \( n \geq m-1 \), then \( f(T)
  \) counts the number of nodes in $T$ and thus gives the \emph{space
  requirement} functional discussed in~\cite[\S3.4]{MR93f:68045}.
\end{example}

\begin{example}
  \label{example:path-length}
  If we define \( f(T) := 0 \) when \( 0 \leq |T| \leq m-2 \) and
  \( t_n := n-(m-1) \) for \( n \geq m-1 \) then \( f \) is the
  \emph{internal path length} functional discussed
  in~\cite[\S3.5]{MR93f:68045}:\ $f(T)$ is the sum of all root-to-key
  distances in $T$.
\end{example}

\begin{example}
  \label{example:shape-functional}
  As described above, each permutation of $[n]$ gives rise to an
  $m$-ary search tree.  Suppose we place the uniform distribution on
  such permutations. This induces a distribution on $m$-ary search
  trees called the \emph{random permutation model}. Denote its
  probability mass function by~$Q$.  It is important to
  note that~$Q$ is \emph{not} uniform, since different permutations
  can give rise to the same tree.    For example, the permutations
  \[(10, 7, 12, 4, 1, 8, 5, 6, 9, 14, 11, 2, 15, 13, 3)\] and
  \[(7, 10,  12, 1, 4, 8, 5, 6, 9, 14, 11, 2, 15, 13,
  3)\] both give rise to the quaternary search tree shown in
  Figure~\ref{fig:terms}.  Dobrow and Fill~\cite{MR97k:68038} noted
  that
  \begin{equation}
    \label{eq:2.2}
    Q(T) = \frac{1}{\prod_{x} \binom{|T_x|}{m-1}},
  \end{equation} where the product in~\eqref{eq:2.2} is over all full
  nodes in \( T \).  This functional is sometimes called the
  \emph{shape functional} as it serves as a crude measure of the
  ``shape'' of the tree, with ``full'' trees (like the complete tree)
  achieving larger values of $Q$.  For further discussions along these
  lines, consult~\cite{MR97k:68038} and~\cite{MR97f:68021}.  If we
  define \( f(T) := 0\) for \( 0 \leq |T| \leq m-2 \) and \( t_n :=
  \ln{\binom{n}{m-1}} \) for \( n \geq m-1 \), then \( f(T) =
  -\ln{Q(T)} \).  In this work we will refer to $-\ln{Q}$ (rather than
  $Q$) as the shape functional.  It was our desire to understand the
  distribution of the shape functional under the random permutation
  model that led to this paper.
\end{example}

\medskip

We study the distribution of the functional \( f(T) \) when $T$ is
given the distribution~$Q$ described in
Example~\ref{example:shape-functional}.  To do this, we derive
asymptotics for the moments of $f(T)$ and then employ the method of
moments.

\begin{remark*}
  In the sequel (without any loss of applicability) we will restrict
  attention to real-valued toll sequences.
\end{remark*}

\subsection*{Related work}
\label{sec:related-work}
Chern and Hwang~\cite{MR1871558} carried out the program of this paper
for the space requirement (Example~\ref{example:space-requirement})
and independently discovered
Theorem~\ref{T:att}~\cite[Proposition~7]{MR1871558}.  Hwang and
Neininger~\cite{01808498} analyzed a range of (not necessarily
deterministic) toll functions for binary search trees.  In this paper
we treat $m$-ary search trees for any $m \geq 2$, and for simplicity
we restrict attention to deterministic toll sequences.  Another
closely related paper is~\cite{MR1933199}, where the authors consider
variants of the \texttt{Quicksort} algorithm by allowing more
general schemes of choosing the pivot.  As in our
Section~\ref{sec:transfer-theorems}, their paper treats
Cauchy--Euler differential equations, and a
generalization of our Theorem~\ref{T:att} is
obtained~\cite[Theorem~1]{MR1933199}.

We obtain moment asymptotics for our additive functionals using the
Asymptotic Transfer Theorem of Section~\ref{sec:transfer-theorems}.
An alternative 
approach employs singularity analysis~\cite{MR90m:05012} of generating
functions.  A sketch of this approach in the case of binary search
trees ($m=2$) is presented in~\cite{FFK}.  One small advantage of the
present approach is that the conditions we impose on the toll sequence
[for example,~\eqref{5.1}] are milder than those required for the
application of singularity analysis.

\section{Overview: main results}
\label{S:intro}

For a given toll sequence, the distribution of $f(T)$ depends only
on~$n$.  We let~$X_n$ denote a random variable whose distribution is that
of $f(T)$ under the random permutation model on~$T$.  It is the
distribution of $X_n$ which is the main focus of this paper.

\subsection{A common framework for all moments}
\label{sec:common-framework-all}
{\sloppy
Under the random permutation model the joint distribution of the
subtree sizes $|L_1|,\dots,|L_m|$ is uniform over all \(
\binom{n}{m-1} \) \( m \)-tuples  of
nonnegative integers that sum to \( n - (m-1)
\):\ see~\cite[Exercise~3.8]{MR93f:68045}.  We now apply the law of
total expectation to compute $\mu_n(k) := \E{X_n^k}$ by conditioning
on $(|L_1|,\ldots,|L_m|)$. 
Let \( \sum_{\mathbf{j}} \) denote the sum over
all \( m \)-tuples \(( j_1,\ldots,j_m) \) that sum to \( n-(m-1) \)
and \( \sum_{\mathbf{k}} \) the sum over all \( (m+1) \)-tuples \( (
k_1,\ldots,k_{m+1} ) \) of nonnegative integers that sum to \( k \).
Then, letting \( \oplus \) denote sums of mutually independent random
variables, for \( n \geq m-1 \) we have
\begin{align*}
  \mu_n(k) &= \E{X_n^k} = \E \E( X_n^k\; \bigl|\; |L_1|,\ldots,|L_m| ) =
  \frac{1}{\binom{n}{m-1}} \sum_{\mathbf{j}} \E{(
  X_{j_1} \oplus \cdots \oplus X_{j_m} + t_n )^k} \\
  &= \frac{1}{\binom{n}{m-1}} \sum_{\mathbf{j}}{\sum_{\mathbf{k}}}
  \binom{k}{k_1,\ldots,k_m,k_{m+1}} \mu_{j_1}(k_1)\cdots
  \mu_{j_m}(k_m) t_n^{k_{m+1}}.
\end{align*}%
}%
We can rewrite this equation as
\begin{equation}
  \label{eq:2.3}
  \mu_n(k) = \frac{m}{\binom{n}{m-1}} \sum_{j=0}^{n-(m-1)} \binom{n-1-j}{m-2}
  \mu_j(k) + r_n(k),
\end{equation}
where
\begin{equation}
  \label{eq:2.4}
  r_n(k) := \sumstark \binom{k}{k_1,\ldots,k_m,k_{m+1}}
  t_n^{k_{m+1}}  \frac{1}{\binom{n}{m-1}} \sum_{\mathbf{j}}
  \mu_{j_1}(k_1) \cdots \mu_{j_m}(k_m),
\end{equation}
with \( \sum_{\mathbf{k}}^{*} \) denoting the same sum as
\(\sum_{\mathbf{k}} \) with the additional restriction that \( k_i < k
\) for \( i=1,\ldots,m \).  We have thus established that the moments
$\mu_n(k)$ each satisfy the same basic recurrence in $n$, differing
as $k$ varies only in the non-homogeneous term \( r_n(k) \).  Observe
that \( r_n(1) = t_n \), the toll function.  We record this important
fact as
\begin{proposition}
  \label{thm:rand-perm-model}
  Under the random permutation model, all moments of an additive
  functional satisfy the basic recurrence
  \begin{equation}
    \label{eq:2.7}
    a_n = b_n + \frac{m}{\binom{n}{m-1}} \sum_{j=0}^{n-(m-1)}
    \binom{n-1-j}{m-2} a_j, \qquad n \geq m-1,
  \end{equation} with specified initial conditions (say) \( a_j := b_j
  \), \( 0 \leq j \leq m-2 \).  [Recall the statement following
  Definition~\ref{def:additive-functional} about the initial
  conditions for the recurrence~\eqref{eq:2.1}.]
\end{proposition}
To be more specific, equation~\eqref{eq:2.7} is satisfied by $a_n =
 \mu_n(k) = \E{X_n^k}$ and $b_n  = r_n(k)$, where $r_n(k)$ is defined
 in terms of lower-order moments of smaller trees at~\eqref{eq:2.4}.
 We proceed to study the recurrence relation~\eqref{eq:2.7} for general input
 $(b_n)$ and corresponding output $(a_n)$.

\subsection{Transfer theorems}
\label{sec:transfer-theorems}

In order to analyze the recurrence relation~\eqref{eq:2.7} we introduce
generating functions
\begin{equation*}
  A(z) := \sum_{n=0}^\infty a_n z^n \quad \text{ and } \quad
  B(z) := \sum_{n=0}^\infty b_n z^n.
\end{equation*}
Furthermore, let
$\rising{x}{r} := \prod_{k=0}^{r-1} (x+k)$
denote the $r$th rising factorial power of $x$ and
$\falling{x}{r} := \prod_{k=0}^{r-1} (x-k)$
the $r$th falling factorial power of $x$.  Multiplying~\eqref{eq:2.7} by \(
\falling{n}{m-1} z^{n-(m-1)} \) and summing over \( n \geq m-1 \) we
get, after some (routine) calculation, the differential equation
\begin{equation}
  \label{eq:2.8}
  A^{(m-1)}(z) = B^{(m-1)}(z) + m! (1-z)^{-(m-1)} A(z).
\end{equation}
Equations of the form~\eqref{eq:2.8} are members of a class known as
\emph{Euler differential equations}.  Using the method of variation of
parameters and combinatorial identities one can obtain the general
solution to this equation.  See Section~\ref{S:ett} for a proof.
\begin{theorem}[Exact Transfer Theorem (ETT)]
  \label{T:ett}
  Let \( A \) and \( B \) denote the
  respective ordinary generating functions of the sequences \( (a_n)
  \) and \( (b_n) \) in the recurrence~\eqref{eq:2.7}.
  Let
\begin{equation}
\label{Bhatdef}
\Bh(z) := B(z) - \sum_{j = 0}^{m - 2} b_j z^j = \sum_{n = m - 1}^{\infty} b_n
z^n.
\end{equation}
Then
\begin{eqnarray}
\label{1.2e}
\hspace{-1in}A(z)
  &=& \sum_{j = 1}^{m - 1} c_j (1 - z)^{- \gl_j} + \sum_{j = 1}^{m - 1} \frac{(1
       - z)^{- \gl_j}}{\psi'(\gl_j)} \int_{0}^z\!B^{(m - 1)}(\zeta) (1 -
       \zeta)^{\gl_j + m - 2}\,d\zeta \\
\label{etta}
  &=& \sum_{j = 1}^{m - 1} c_j (1 - z)^{- \gl_j} + \Bh(z) + m! \sum_{j = 1}^{m -
       1} \frac{(1 - z)^{- \gl_j}}{\psi'(\gl_j)} \int_{0}^z 
\Bh(\zeta) (1
       - \zeta)^{\gl_j - 1}\,d\zeta,
\end{eqnarray}
  where \( \psi \) is the indicial polynomial
  \begin{equation}
    \label{eq:2.10}
    \psi(\lambda) := \rising{\lambda}{m-1} - m! = \lambda(\lambda+1)
    \cdots (\lambda+m-2) - m!
  \end{equation}
  with roots \( 2 =: \lambda_1,\lambda_2, \ldots,
  \lambda_{m-1} \) in nonincreasing order of real parts.
  In~\eqref{1.2e}, the
  coefficients  \( c_1, c_2, \ldots, c_{m-1} \)  can
  be written explicitly in terms of the initial conditions \( b_0,
  \ldots, b_{m-2} \) as
  \begin{equation}
    \label{eq:2.9}
    c_j = \frac{m!}{\psi'(\lambda_j)} \sum_{k=0}^{m-2} b_k
    \frac{k!}{\rising{\lambda_j}{k+1}}, \qquad j=1,\ldots,m-1.
  \end{equation}
\end{theorem}
\noindent
In particular,
\begin{equation}
\label{2.1}
c_1 = \frac{1}{H_m - 1} \sum_{j = 0}^{m - 2} \frac{b_j}{(j + 1) (j + 2)}.
\end{equation}

The indicial polynomial~\eqref{eq:2.10} is well studied;
see~\cite{MR90a:68012,MR93f:68045,MR1871558} and
Appendix~\ref{sec:prop-indic-polyn}, and also the
related~\cite{MR96j:68042}.  
We will exploit the
expression~\eqref{1.2e} for \( A(z) \) to relate the asymptotic
properties of the sequence \( (b_n) \) to those of \( (a_n) \) and
then use \emph{transfer theorems} to derive limiting distributions of
additive functionals.


\begin{remark}
For computations, equation~\eqref{1.2e} might be easier to use when it is no
bother to compute derivatives of~$B$; otherwise, \eqref{etta} is easier.
Equation~\eqref{etta} will be used in establishing part~(a) of the Asymptotic
Transfer \refT{T:att}; the proof of part~(b) will use~\eqref{1.2e}.
\end{remark}

\medskip
It is quite easy to transfer asymptotics for~$B$ to asymptotics for~$A$ using
the ETT.  We give three examples important for applications to moments of
functionals in the next theorem, proved independently by F\'{e}lix Chern and
Hsien-Kuei Hwang in~\cite{MR1871558} using a quite
different approach.  The series convergence required in~\eqref{5.1}
need not be absolute.

\begin{theorem}[Asymptotic Transfer Theorem (ATT)]
\label{T:att}
\
\smallskip

{\rm (a)}~If
\begin{equation}
\label{5.1}
b_n = o(n) \qquad \mbox{\rm and} \qquad \sum_{n = 0}^{\infty} 
\frac{b_n}{(n + 1)
(n + 2)}\mbox{\rm \ converges,}
\end{equation}
then
\begin{equation}
\label{5.1a}
a_n = \frac{K_1}{H_m - 1} n + o(n), \qquad \mbox{\rm where} \qquad K_1 :=
\sum_{j = 0}^{\infty} \frac{b_j}{(j + 1) (j + 2)}.
\end{equation}
\smallskip

{\rm (b)}~If
$b_n \equiv K_2 (n + 1) + h_n$ where $(h_n)$ satisfies~\eqref{5.1}
[with $(b_n)$ replaced by $(h_n)$], then
\begin{equation}
\label{5.2}
a_n = \frac{K_2}{H_m - 1} n H_n + \frac{K_3}{H_m - 1} n + o(n),
\end{equation}
where
\begin{equation}
\label{5.3}
K_3 := \sum_{j = 0}^{\infty} \frac{h_j}{(j + 1) (j + 2)} + K_2 \left[ \frac{H_m
- 1}{2} - 1 + \frac{H^{(2)}_m - 1}{2 (H_m - 1)} \right].
\end{equation}
\smallskip

{\rm (c)}~If $b_n = K_4 n^v + o\left(n^{v}\right)$ with $v > 1$,
then
\begin{equation}
\label{5.4}
a_n = \frac{K_4}{1 - \frac{m! \Gamma(v + 1)}{\Gamma(v + m)}} n^v +
o\left(n^{v}\right).
\end{equation}
\end{theorem}

Of course part~(a) is equally valid with $\sum_{n = 0}^{\infty}
\frac{b_n}{n^2}$ replacing the series in~\eqref{5.1}.

\refT{T:att} will be proved in \refS{S:att}.  It is easy to see that the ratio
$\frac{m! \Gamma(v + 1)}{\Gamma(v + m)}$ appearing in~\eqref{5.4} has modulus
strictly less than unity; in particular, the expression given is well defined.
Refined and additional asymptotic transfers will be discussed in 
Section~5.


\subsection{Limiting distributions}
\label{sec:limit-distr}

Applications to moments and
limiting distributions are discussed for ``small'' toll sequences in
\refS{S:small}---see Theorems~\ref{T:small1} and~\ref{T:small2};
``moderate'' and ``large'' toll sequences are discussed in
Section~\ref{sec:moderate-large-toll}---see Theorems~\ref{thm:F2}
and~\ref{thm:F3}.
We state here a summary theorem that can be deduced easily from these
four theorems.  [For the definition of $m_0(\beta)$ in case~\ref{item:2} of
the theorem, see Remark~\ref{remark:alpha-monotone}.  Concerning the
values of $g_1$ and $g_2$ in cases~\ref{item:2} and~\ref{item:3},
consult~\eqref{eq:8}, \eqref{eq:9} (and~\eqref{eq:2}).]
\begin{theorem}[Limit theorem for additive functionals]
  \label{thm:limit-summary}
  Let $X_n = f(T_n)$ be the additive functional on random $m$-ary
  search trees corresponding to a toll sequence
  $(t_n)$, with $T_n$ having the distribution on trees induced by a uniformly
  distributed permutation on $[n]$.  When $(t_n)$ satisfies
  one of the three
  conditions~\ref{item:1}--\ref{item:3} below, then
  \begin{equation*}
    \frac{X_n - \E{X_n}}{\sqrt{\Var{X_n}}} \xrightarrow{\mathcal{L}}
    W, 
  \end{equation*}
  with convergence of all moments.  Here, with $L$ denoting a
  slowly varying function,
  \begin{enumerate}
  \item If $2 \leq m \leq 26$ and
    {\renewcommand{\labelenumii}{{\normalfont (\alph{enumii})}}%
    \begin{enumerate}
    \item $t_n = o(\sqrt{n})$ and $\sum^{\infty} n^{-1} \max_{n^\delta
        \leq k \leq n} \frac{t_k^2}{k^2} < \infty$ for some $0 <
      \delta < 1$ or
    \item $t_n \sim \sqrt{n}L(n)$,
    \end{enumerate}%
    }
    then $W$ has the standard normal distribution.
    \label{item:1}
  \item     \label{item:2}
    If $t_n \sim n^\beta L(n)$, with $1/2 < \beta < 1$ and $2 \leq m
    \leq m_0(\beta)$, then $W = g_2^{-1/2}Y$ where $\mathcal{L}(Y)$ is the
    unique distribution with finite second moment
    satisfying~\eqref{eq:1}.  Here
    \begin{equation}
      \label{eq:14}
      g_2 = \frac{(m-1)!}{1 - \frac{m! \Gamma(2\beta+1)}{\Gamma(2\beta
      + m)}} \left[ \frac{1}{(m-1)!} + 2m g_1
      \frac{\Gamma(\beta+1)}{\Gamma(\beta + m)} + m(m-1) g_1^2
      \frac{\Gamma^2(\beta+1)}{\Gamma(2\beta + m)} \right] > 0,
    \end{equation}
    with
    \begin{equation}
      \label{eq:15}
      g_1 =
      \left(
        1 - \frac{m! \Gamma(\beta+1)}{\Gamma(\beta+m)}
      \right)^{-1}.
    \end{equation}
  \item If $t_n \sim n^\beta L(n)$ with $\beta > 1$, then $W =
    g_2^{-1/2}Y$ where $\mathcal{L}(Y)$
    is the unique distribution with finite second moment
    satisfying~\eqref{eq:1}, where $g_2$ is again defined
    at~\eqref{eq:14}--\eqref{eq:15}.
    \label{item:3}
  \end{enumerate}
\end{theorem}
\begin{remark}
  \label{remark:g_1-g_2}
  In case~\ref{item:3}, it is easy to check that $g_1 > 0$ and $g_2 >
  0$.  In case~\ref{item:2} one can verify easily that $g_1 < 0$.  In
  this case to see that $g_2 > 0$, note that no constant random
  variable satisfies~\eqref{eq:1}.
\end{remark}

\section{Proof of the ETT}
\label{S:ett}

In this section we prove the ETT, which is \refT{T:ett}.
\begin{proof}
  For the proof of~\eqref{1.2e}, see Appendix~\ref{appendix:diffeq},
  in  particular 
  Corollary~\ref{corollary:A.5.2} and Proposition~\ref{lemma:A.2.1}.
To begin the proof of~\eqref{etta}, note that~$B$
can be replaced by~$\Bh$ in~\eqref{1.2e}.  We then use repeated
  integration by parts and 
invoke Identity~\ref{identity:B.2.2}.  Denoting
\[
\Ah := A(z) - \sum_{j = 1}^{m - 1} c_j (1 - z)^{- \gl_j},
\]
after $m - 2$ integrations by parts we find
\[
\Ah(z) = \sum_{j = 1}^{m - 1} \frac{(1 - z)^{- \gl_j}}{\psi'(\gl_j)} (\gl_j + m
- 2) \cdots (\gl_j + 1) \int_{0}^z\!\Bh'(\zeta) (1 -
\zeta)^{\gl_j}\,d\zeta.
\]
But
\[
(\gl_j + m - 2) \cdots (\gl_j + 1) = \frac{\rising{\gl_j}{m -
1}}{\gl_j} = \frac{\psi(\gl_j) + m!}{\gl_j} = \frac{m!}{\gl_j},
\]
so
\[
\Ah(z) = m! \sum_{j = 1}^{m - 1} \frac{(1 - z)^{- \gl_j}}{\gl_j \psi'(\gl_j)}
\int_{0}^z\!\Bh'(\zeta) (1 - \zeta)^{\gl_j}\,d\zeta.
\]
We obtain~\eqref{etta} by performing one more integration by parts and
utilizing Identity~\ref{identity:B.1} with $\gl = 0$.
%
\end{proof}

\section{Proof of the ATT}
\label{S:att}
In this section we prove the ATT, which is \refT{T:att}.  The  following
elementary result, which is found in the first line of the proof of
Lemma~6 in~\cite{MR1871558}, is key to the analysis of~\eqref{etta}.  For
completeness, we include a proof.

\begin{lemma}
\label{L:exactlinear}
Let $Y(z) = \sum_{n = 0}^{\infty} y_n z^n$ with $y_0 = 0$.  For any $\gl \in
\CC$,
\begin{equation}
\label{exactlinear}
[z^n] \left( (1 - z)^{-\gl} \int_{0}^z\!(1 - \zeta)^{\gl - 1}
Y(\zeta)\,d\zeta \right) = \sum_{k = 0}^{n - 1} \frac{y_k}{k + 1} 
\prod_{j = k +
2}^n \left( 1 + \frac{\gl - 1}{j} \right), \qquad n \geq 0.
\end{equation}
\end{lemma}
\noindent
The product in~\eqref{exactlinear} may be written (when $\gl \in \CC \setminus
\{0, -1, -2, \ldots\}$)  as
\[
\frac{\Gamma(\gl + n) \Gamma(2 + k)}{\Gamma(1 + n) \Gamma(\gl + k + 1)},
\]
which by Stirling's formula equals
\begin{equation}
\label{Stirling}
\frac{n^{\gl - 1} \left[ 1 + O(n^{-1}) \right]}{(k + 1)^{\gl - 1} \left[ 1 +
O((k + 1)^{-1}) \right]}
\end{equation}
for $n \geq 1$ and $k \geq 1$.  [The product in~\eqref{exactlinear}
equals~\eqref{Stirling} even if $\gl \in \{0, -1, -2, \ldots\}$.]  Also, of
special interest is the case $\gl = 2$, in which case~\eqref{exactlinear}
reduces to
\begin{equation}
\label{exactlinear2}
[z^n] \left( (1 - z)^{-2} \int_{0}^z\!(1 - \zeta) Y(\zeta)\,d\zeta
\right)  = (n + 1) \sum_{k = 0}^{n - 1} \frac{y_k}{(k + 1) (k + 2)}, \qquad n
\geq 0.
\end{equation}

\begin{proof}
The function $W(z) := (1 - z)^{-\gl} \int_{0}^z\!(1 - \zeta)^{\gl -
1} Y(\zeta)\,d\zeta$ is the unique solution with $W(0) = 0$ to the differential
equation
\[
W'(z) = \gl (1 - z)^{-1} W(z) + (1 - z)^{-1} Y(z);
\]
that is, $w_n := [z^n] W(z)$, $n \geq 0$, satisfies $w_0 = 0$ and
\begin{equation}
\label{wn}
w_n = \frac{\gl}{n} \sum_{k = 0}^{n - 1} w_k + \frac{1}{n} \sum_{k = 0}^{n - 1}
y_k, \qquad n \geq 1.
\end{equation}
But the recurrence~\eqref{wn} can be easily solved to
yield~\eqref{exactlinear}:\ compute the difference $n w_n - (n - 1)
w_{n - 1}$ and iterate.
\end{proof}

For part~(a) of the ATT we use the following estimates from~\cite{MR1871558};
these follow readily from \eqref{exactlinear}--\eqref{exactlinear2}.  [Part~(a)
is their Lemma~6; part~(b) is used tacitly in the proof of their 
Proposition~7.]

\begin{lemma}
\label{L:linear}
\
\smallskip

{\rm (a)}~If $\Re(\gl) < 2$ and $Y(z) = \sum_{n = 0}^{\infty} y_n z^n$
satisfies $y_0 = 0$ and $y_n = o(n)$, then
\[
[z^n] \left( (1 - z)^{- \gl} \int_{0}^z\!Y(\zeta) (1 - \zeta)^{\gl -
1}\,d\zeta \right) = o(n).
\]
\smallskip

{\rm (b)}~With~$\Bh$ defined at~\eqref{Bhatdef}, if~\eqref{5.1} holds, then
\[
[z^n] \left( (1 - z)^{-2} \int_{0}^z\!\Bh(\zeta) (1 - \zeta)\,d\zeta
\right) = n \sum_{j = m - 1}^{\infty} \frac{b_j}{(j + 1) (j + 2)} + o(n).~\nopf
\]
\end{lemma}

For part~(c) of the ATT, we will
need the following simple comparison lemma.

\begin{lemma}
\label{L:comp}
If $(b_n)$ and $(b'_n)$ are two input sequences such that
\[
|b_n| \leq b'_n\mbox{\rm\ \ for all~$n \geq 0$},
\]
then the corresponding output sequences $(a_n)$ and $(a'_n)$ in~\eqref{eq:2.7}
(with the initial conditions stated there) satisfy
\[
|a_n| \leq a'_n\mbox{\rm\ \ for all~$n \geq 0$}.
\]
\end{lemma}

\begin{proof}
This follows immediately by induction.
\end{proof}

\begin{proof}[Proof of \refT{T:att}]
(a)~From~\eqref{etta}, assumption~\eqref{5.1}, \refL{L:linear}, 
\eqref{2.1}, and
$\psi'(2) = m! (H_m - 1)$, the result is immediate.

(b)~Suppose first that $b_n \equiv n + 1$.  Then $B(z) \equiv (1 - z)^{-2}$,
so
\[
B^{(m - 1)}(z) \equiv m! (1 - z)^{- (m + 1)}.
\]
Plugging this into~\eqref{1.2e} we find
\[
a_n = (n + 1) \left[ c_1 + m! \sum_{j = 2}^{m - 1} \frac{1}{(2 - \gl_j)
\psi'(\gl_j)} \right] + \frac{m!}{\psi'(2)} [z^n] \left[ (1 - z)^{-2}
\log\left((1 - z)^{-1} \right) \right] + o(n).
\]
Now we use Identities~\ref{identity:psp2} and~\ref{identity:B.5.22}.
Also, since $b_n \equiv n + 1$, we have from~\eqref{2.1} in this case
that $c_1 = 1$.  Therefore
\begin{eqnarray*}
a_n
  &=& (n + 1) \left( 1 + \frac{1}{2} \left[ \frac{H^{(2)}_m - 1}{(H_m - 1)^2} -
       1\right] \right) + \frac{1}{H_m - 1} [(n + 1) H_n - n] + o(n) \\
  &=& \frac{1}{H_m - 1} n H_n + \left[ \frac{1}{2} - \frac{1}{H_m - 1} +
\frac{H^{(2)}_m - 1}{2 (H_m - 1)^2} \right] n + o(n).
\end{eqnarray*}
This completes the proof of~(b) for our special case, and the general case
follows from this and part~(a) using the superposition principle.

(c)~Suppose first that $b_n \equiv \rising{(v + 1)}n / n! \sim n^v /
\Gamma(v + 1)$, so that
$B(z)
\equiv (1 - z)^{-(v + 1)}$ and $B^{(m - 1)}(z) \equiv \rising{(v + 1)}{m
- 1} (1 - z)^{-(v + m)}$.  Plugging this into~\eqref{1.2e} and
utilizing Identity~\ref{identity:B.1} with $\gl = v + 1$ and the calculation
\[
\frac{\rising{(v + 1)}{m - 1}}{\rising{(v + 1)}{m - 1} - m!} =
\left[ 1 - \frac{m! \Gamma(v + 1)}{\Gamma(v + m)} \right]^{-1},
\]
we find
\[
A(z) = \left[ 1 - \frac{m! \Gamma(v + 1)}{\Gamma(v + m)} \right]^{-1} (1 -
z)^{- (v + 1)} + O(|1 - z|^{-2}).
\]
By singularity analysis for large functions~\cite{MR90m:05012}, this
completes the proof of~(c) for our special case.

To complete the proof in the general case, we need only show that if $b_n =
o(n^v)$ for $v > 1$, then $a_n = o(n^v)$.  Indeed, fix $\eps > 0$; then
there exists a sequence $(b'_n)$ such that $|b_n| \leq b'_n$ for all~$n$
and
\[
b'_n = \eps \rising{(v + 1)}n / n!\mbox{\ \ for all large~$n$}.
\]
The toll sequence is but a slight modification of our special-case toll
sequence, and we see that
\[
a'_n = \eps' n^v + o(n^v),\mbox{\ \ where\ \ }\eps' := \frac{\eps}{\Gamma(v +
1)} \left[ 1 - \frac{m! \Gamma(v + 1)}{\Gamma(v + m)} \right]^{-1}.
\]
Now \refL{L:comp} implies that
\[
\limsup_n |a_n| n^{-v} \leq \eps';
\]
since~$\eps$ (and hence $\eps'$) can be made arbitrarily small, this completes
the proof.
\end{proof}

The conditions~\eqref{5.1} on $(b_n)$ are not only sufficient but 
also necessary
for asymptotic linearity of~$a_n$.  Indeed, here is a converse:

\begin{proposition}
\label{converse}
If $a_n = K n + o(n)$ for some constant~$K$, then~\eqref{5.1} holds.
\end{proposition}

\begin{proof}
  Chern and Hwang~\cite{MR1871558} provide the simple proof that $b_n
  = o(n)$.  Moreover, then, from~\eqref{etta}, \eqref{2.1},
  and~\eqref{exactlinear2} with~$Y$ taken to be~$\Bh$,
  we find that
\[
\sum_{j = 0}^{n - 1} \frac{b_j}{(j + 1) (j + 2)} = K (H_m - 1) + o(1),
\]
i.e.,\ that the series $\sum_{n = 0}^{\infty} \frac{b_n}{(n + 1) (n +
  2)}$ converges [to $K (H_m - 1)$].
\end{proof}

The following additional asymptotic transfer results are established by
calculations similar to those in the proof of the ATT.  We leave detailed
proofs as exercises for the reader.

\begin{theorem}[more asymptotic transfers]
\label{T:more}
Consider the initial value problem~\eqref{eq:2.7}.
\smallskip

{\rm (a)}~If $2 \leq m \leq 26$ and $b_n = o(\sqrt{n})$, then we can
refine~\eqref{5.1a} to
\begin{equation}
\label{smallremainder}
a_n = \frac{K_1}{H_m - 1} n + o\left(\sqrt{n}\right).
\end{equation}
\smallskip

{\rm (b)}~If
$\Re\,\lambda_2 =: \alpha < 1 + \beta$ and $b_n \sim n^\beta L(n)$,
where $1/2 < \beta < 1$,
with~$L$ slowly
varying, then we can refine~\eqref{5.1a} to
\begin{equation}
\label{slowremainder}
a_n = \frac{K_1}{H_m - 1} n - \frac{\rising{(1+\beta)}{m -
    1}}{m! - \rising{(1+\beta)}{m - 1}} n^\beta L(n) + o\left(n^\beta
    L(n)\right). 
\end{equation}

{\rm (c)}~If $b_n \sim n L(n)$ with~$L$ slowly varying, then, with $K_1$
defined at~\eqref{5.1a},
\begin{equation}
\label{linearslowlead}
a_n \sim
  \left\{ \begin{array}{cl}
    \frac{K_1}{H_m - 1} n, & \mbox{\rm if\ }\sum^{\infty} \frac{L(k)}{k} <
                               \infty, \\
    \frac{1}{H_m - 1} n \sum_{k \leq n} \frac{L(k)}{k},
                           & \mbox{\rm if\ }\sum^{\infty} \frac{L(k)}{k} =
                               \infty.
           \end{array}
   \right.
\end{equation}

{\rm (d)}~Part~{\rm (c)} of the {\rm ATT} can be extended as follows.  If $b_n
= K_4 n^v L(n) + o\left(n^{v} L(n)\right)$ with $v > 1$ and~$L$
slowly varying, then
\begin{equation}
\label{higherslowlead}
a_n = \frac{K_4}{1 - \frac{m! \Gamma(v + 1)}{\Gamma(v + m)}} n^v L(n) +
o\left(n^{v} L(n)\right).
\end{equation}

\end{theorem}

\begin{proof}[Proof hints]
Whenever the conditions~\eqref{5.1} are met we have by the ETT
and~\eqref{exactlinear2}
\begin{eqnarray}
\label{muntilde}
a_n - \frac{K_1}{H_m - 1} (n + 1)
  &=& O(n^{\alpha - 1}) + b_n - \frac{1}{H_m - 1} (n + 1) \sum_{k = n}^{\infty}
        \frac{\bh_k}{(k + 1) (k + 2)} \\
  & & \hspace{.1in} + m! \sum_{j = 2}^{m - 1} \frac{1}{\psi'(\gl_j)} 
[z^n] \left(
        (1 - z)^{- \gl_j} \int_{0}^z\!\Bh(\zeta) (1 - \zeta)^{\gl_j -
        1}\,d\zeta \right), \nonumber
\end{eqnarray}
where $\alpha$ is strictly smaller than $1 + \beta$ by
  assumption.
(When $\beta = 1/2$ we know that $\alpha < 3/2$ when $m \leq 26$.)
Simple estimates, including the use of~\eqref{Stirling},
give cases~(a) and~(b); for~(b), the coefficient of $n^\beta L(n)$
in~\eqref{slowremainder} indeed is, using
  Identities~\ref{identity:B.1} and~\ref{identity:psp2},
\[
1 - \frac{1}{(1 - \beta)(H_m - 1)} + m! \sum_{j = 2}^{m - 1}
\frac{1}{((1+\beta) - \gl_j)
\psi'(\gl_j)} = - \frac{\rising{(1+\beta)}{m - 1}}{m! - \rising{(1 +
  \beta)}{m - 1}}. 
\]

In case~(c), from the ETT result~\eqref{etta} and simple estimates we find
\[
a_n = (n + 1) \frac{1}{H_m - 1} \sum_{k = 0}^{n - 1} \frac{b_k}{(k + 1) (k +
2)} + O(n L(n)).
\]
To finish, we use the regular variation fact (quoted by Hwang and
Neininger~\cite{01808498}
at their equation~(7) from Proposition 1.5.9a in~\cite{MR90i:26003}) that
\begin{equation}
\label{bingham}
L(n) = o\left(\sum_{k \leq n} \frac{L(k)}{k}\right).
\end{equation}

In case~(d), again from~\eqref{etta} and simple estimates we find
\[
a_n = O(n) + b_n + K_4 n^v L(n) m! \sum_{j = 1}^{m - 1}
\frac{1}{\psi'(\gl_j) (v + 1 - \gl_j)} + o\left(n^v L(n)\right).
\]
The proof is completed by using Identity~\ref{identity:B.1}.
\end{proof}

\section{Asymptotic normality for small toll functions}
\label{S:small}

In this section we establish asymptotic normality for ``small'' toll
functions.

\subsection{Central limit theorem statements}
\label{S:CLT}

As did Hwang and Neininger in~\cite{01808498}, in our Theorems~\ref{T:small1}
and~\ref{T:small2} we treat two overlapping cases.  Throughout, we write
$\sum_{\jj}$ as shorthand for the sum over $m$-tuples $(j_1, \ldots, j_m)$ of
nonnegative integers summing to $n - (m - 1)$.

\begin{theorem}[CLT I for small toll functions]
\label{T:small1}
If 
$2 \leq m \leq 26$ and the toll sequence $(t_n)$ satisfies
\begin{equation}
\label{small1}
\mbox{\rm (a)\ }t_n = o\left(\sqrt{n}\right) \qquad \mbox{\rm and} \qquad
\mbox{\rm (b)\ }\sum^{\infty} n^{-1} \max_{n^{\gd} \leq k \leq n}
\frac{t^2_k}{k} < \infty\mbox{\rm \ \ for some $0 < \gd < 1$},
\end{equation}
then the mean $\mu_n$ and variance $\sigma^2_n$ of the corresponding
additive functional~$X_n$ on $m$-ary search trees with the random
permutation model satisfy, respectively,
\begin{equation}
\label{asymean}
\mu_n = \frac{K_1}{H_m - 1} n + o\left(\sqrt{n}\right) =: \mu n +
o\left(\sqrt{n}\right),
\end{equation}
with~$K_1$ defined at~\eqref{5.1a}, and
\begin{equation}
\label{asyvar}
\sigma^2_n = \sigma^2 n + o(n), \qquad \mbox{\rm where} \qquad \sigma^2 :=
\frac{1}{H_m - 1} \sum_{j = 0}^\infty \frac{r_j}{(j + 1) (j + 2)},
\end{equation}
with the sequence $(r_n)$ defined by $r_j := 0$ for $0 \leq j \leq
m - 2$ and
\begin{equation}
\label{rdef}
r_n := \frac{1}{{n \choose {m - 1}}} \sum_{\jj} \big[ t_n + \mu_{j_1} +
\cdots + \mu_{j_m} - \mu_n \big]^2, \qquad n \geq m - 1.
\end{equation}
Moreover,
\[
\frac{X_n - \mu n}{\sqrt{n}} \mbox{\rm \ \ is asymptotically\ \ N$(0,
\sigma^2)$},
\]
and there is convergence of moments of every order.
\end{theorem}

\begin{remark}
\label{R:degenerate}
One can check 
(for \emph{any} $2 \leq m < \infty$)
that the variance $\sigma^2_n$ vanishes for all $n \geq m - 1$
if and only if the toll sequence is chosen as
\[
t_n = t \min\{m - 1, n\}, \qquad n \geq 0
\]
for some constant $t \in \RR$.  In that case without loss of generality $t =
1$ and then $X_n \equiv n$ is just the number of keys and we have exact
(though degenerate) normality.  So we shall assume in the proof of
\refT{T:small1} that $\sigma^2 > 0$.
\end{remark}

\begin{remark}
(a)~Condition~\eqref{small1}(b) trivially implies
\begin{equation}
\label{bsimpler}
\sum^{\infty} \frac{t^2_n}{n^2} < \infty,
\end{equation}
which in turn implies that~\eqref{5.1} holds with absolute convergence; indeed,
since the nonnegative numbers $[(n + 1) (n + 2)]^{-1}$, $n \geq 0$, sum to
unity, we have
\begin{equation}
\label{bcs}
\left[ \sum_{n = 0}^{\infty} \frac{|t_n|}{(n + 1) (n + 2)} \right]^2 \leq
\sum_{n = 0}^{\infty} \frac{t^2_n}{(n + 1) (n + 2)} < \infty.
\end{equation}

(b)~If
\begin{equation}
\label{bt}
|t_n| =
O\left(\tilde{t}_n\right)\mbox{\ \ with\ \ }0 \leq
\frac{\tilde{t}_n}{\sqrt{n}} \downarrow\mbox{\ \ and \ \
}\sum^{\infty} \frac{\tilde{t}^2_n}{n^2} < \infty, 
\end{equation}
then we claim that~\eqref{small1} holds, and then as a corollary
\[
\frac{t_n}{\sqrt{n}} \downarrow 0\mbox{\ \ and \ \ }\sum^{\infty}
\frac{t^2_n}{n^2} < \infty
\]
implies~\eqref{small1}.  To see the claim, first observe that the
condition~\eqref{bt} certainly implies~\eqref{small1}(a); moreover, we observe
that the series (say, over $2 \leq n < \infty$) in~\eqref{small1}(b) is bounded
by a constant times
\[
\sum_{n = 2}^{\infty} n^{-1} \max_{n^{\delta} \leq k \leq n}
\frac{\tilde{t}^2_k}{k} 
= \sum_{n = 2}^{\infty} n^{-1} \frac{\tilde{t}^2_{\lc n^{\gd}
    \rc}}{\lc n^{\gd} \rc} 
= \sum_{k = 2}^{\infty} \frac{\tilde{t}^2_k}{k} \sum_{(k - 1)^{1 /
    \gd} < n \leq 
k^{1 / \gd}}\!\!n^{-1} = O\left(\sum_{k = 2}^{\infty} 
\frac{\tilde{t}^2_k}{k^2}\right)
< \infty.
\]

(c)~One can check that, when $m = 2$, the proof we will give of CLT~I requires
only~\eqref{small1}(a) and~\eqref{bsimpler}.  In that case we
obtain a strengthening of ``Case S1'' of Theorem~2 in~\cite{01808498} (for
deterministic toll sequences); they required $t_n = O\left(\sqrt{n} / (\log
n)^{(1 / 2) + \eps}\right)$ for some $\eps > 0$.
\end{remark}

\begin{theorem}[CLT II for small toll functions]
\label{T:small2}
If $2 \leq m \leq 26$ and the toll sequence $(t_n)$ satisfies
\begin{equation}
\label{small2}
t_n \sim \sqrt{n} L(n)
\end{equation}
with~$L$ slowly varying,
then the mean $\mu_n$ of the corresponding additive
functional~$X_n$ on $m$-ary search trees with the random permutation model
satisfies
\begin{equation}
\label{asymean2}
\mu_n = \frac{K_1}{H_m - 1} n - \frac{\rising{(3 / 2)}{m - 1}}{m! -
\rising{(3 / 2)}{m - 1}} \sqrt{n} L(n) + o\left(\sqrt{n} L(n)\right).
\end{equation}
with~$K_1$ defined at~\eqref{5.1a}.
If $\sum^{\infty} \frac{L^2(k)}{k} < \infty$, then the variance~$\sigma^2_n$
satisfies~\eqref{asyvar}--\eqref{rdef} and we define
\[
s^2(n) := \sigma^2 n.
\]
If $\sum^{\infty} \frac{L^2(k)}{k} = \infty$, then
\begin{equation}
\label{asyvar2}
\sigma^2_n \sim s^2(n) := \sigma^2 n \sum_{k \leq n} \frac{L^2(k)}{k}
\end{equation}
where in this case we define
\begin{equation}
\label{sigmadef2}
\sigma^2 := \frac{\left(\rising{(\frac{3}{2})}{m - 1} \right)^2 \left[
\frac{\pi}{4} (m - 1) + 1 \right] - (m!)^2}{(H_m - 1) \left[ m! -
\rising{(\frac{3}{2})}{m - 1} \right]^2}. 
\end{equation}
Moreover, in either case
\[
\frac{X_n - \mu n}{s(n)} \mbox{\rm \ is asymptotically standard normal},
\]
and there is convergence of moments of every order.
\end{theorem}

\begin{remark}
When $m = 2$ the constant $\sigma^2$ in~\eqref{sigmadef2} equals $\frac{9}{2}
\pi - 14$, and \refT{T:small2} reduces to ``Case S2'' of Theorem~2 in~\cite{01808498} (for deterministic toll sequences): see especially their
displays~(15) and~(17), with $\tau_2 = 1$.
\end{remark}

\subsection{Central limit theorem proofs}
\label{S:CLTpfs}

\begin{proof}[Proof of CLT~I (\refT{T:small1})]
We use the method of moments together with asymptotic transfer results.

Given the toll sequence $(t_n)$ defining the sequence $(X_n)$ of
random functionals of interest, the means $(\mu_n)$ satisfy the basic
recurrence relation~\eqref{eq:2.7} with $(b_n)$ replaced by $(t_n)$.
Thus~\eqref{asymean} simply repeats the asymptotic transfer
result~\eqref{smallremainder}.

According to the law of total variance, the sequence $(\sigma^2_n)$ of
variances also satisfies the recurrence~\eqref{eq:2.7}, but with $(b_n)$ replaced
by $(r_n)$.  According to \refL{L:rconds} to follow, the sequence $(r_n)$
satisfies the conditions~\eqref{5.1}.  [Note:\ When $m = 2$,
only~\eqref{small1}(a) and~\eqref{bsimpler}, not~\eqref{small1}(b), are needed
in the proof of \refL{L:rconds}.]  Then~\eqref{asyvar} is immediate from
part~(a) of the ATT.

Let $\Xt_n := X_n - \mu (n + 1)$ for $n \geq 0$.  We will complete the
proof of CLT~I by showing by induction on~$k$ that
\begin{equation}
\label{10.1}
\mut_n(k) := \EE\,\Xt^k_n, \qquad k \geq 1 \qquad \mbox{[with 
$\mut_n(0) := 1$]}
\end{equation}
satisfies
\begin{equation}
\label{even}
\mut_n(2 k) \sim \frac{(2 k)!}{2^k k!} \sigma^{2 k} n^k, \qquad k \geq 1
\end{equation}
and
\begin{equation}
\label{odd}
\mut_n(2 k - 1) = o\left(n^{k - (1 / 2)}\right), \qquad k \geq 1.
\end{equation}
Observe that \eqref{asymean}--\eqref{asyvar} imply that~\eqref{even}
and~\eqref{odd} both hold for $k = 1$.

The key to the induction step for both~\eqref{even} and~\eqref{odd} is
to apply the law of total expectation to~\eqref{10.1}, by conditioning
on the subtree sizes $|L_1|, \ldots, |L_m|$ (recall the notation in
Section~\ref{sec:background}).  In a manner 
similar to \eqref{eq:2.3}, we have
\begin{equation}\label{induction}
\mut_n(k)
  = \frac{m}{{n \choose {m - 1}}} \sum_{j = 0}^{n - (m - 1)} {{n - 1 - j}
        \choose {m - 2}} \mut_j(k) + r_n(k),
\end{equation}
where
\begin{equation}
  \label{eq:11}
r_n(k) := \sideset{}{^*}\sum_{\kk} {k \choose {k_1, \ldots, k_m, k_{m
      + 1}}} t^{k_{m + 1}}_n \times \frac{1}{{n \choose {m - 1}}}
      \sum_{\jj} \mut_{j_1}(k_1) \cdots 
\mut_{j_m}(k_m)
\end{equation}
with $\sum^*_{\kk}$ here denoting the same sum as $\sum_{\kk}$ with
the additional restriction that $k_i < k$ for $i = 1, \ldots, m$.  Observe
that~\eqref{induction} is again of the basic form~\eqref{eq:2.7}.  We will apply
the ATT after evaluating $r_n(k)$ asymptotically.

We will treat the induction step in detail only for~\eqref{even}, the
case~\eqref{odd} being similar and somewhat easier.  Let $\sum^{**}_{\kk}$
denote the sum over $m$-tuples $(k_1, \ldots, k_m)$ of nonnegative integers,
each $< k$, summing to~$k$ (i.e.,\ the same sum as $\sum^*_{\kk}$ with the
additional restriction that $k_{m + 1} = 0$).  For $k \geq 2$ we clearly have,
by induction [recall \eqref{even}--\eqref{odd}],
\begin{eqnarray*}
r_n(2 k)
  &=& o(n^k) + \sideset{}{_{}^{**}}\sum_{\kk} {{2k} \choose {2 k_1, \ldots,
  2 k_m}}
        \frac{1}{{n \choose {m - 1}}} \sum_{\jj} \mut_{j_1}(2 k_1) \cdots
        \mut_{j_m}(2 k_m) \\
  &=& o(n^k) + \sideset{}{_{}^{**}}\sum_{\kk} {{2k} \choose {2 k_1,
  \ldots, 2 k_m}}
        \frac{1}{{n \choose {m - 1}}} \sum_{\jj} \frac{(2 k_1)!}{2^{k_1} k_1!}
        \sigma^{2 k_1} j_1^{k_1} \cdots \frac{(2 k_m)!}{2^{k_m} k_m!}
        \sigma^{2 k_m} j_m^{k_m} \\
  &=& o(n^k) + \frac{(2 k)!}{2^k k!} \sigma^{2 k} n^k
  \sideset{}{_{}^{**}}\sum_{\kk} {k \choose
        {k_1, \ldots, k_m}} \frac{1}{{n \choose {m - 1}}} \sum_{\jj}
        \left( \frac{j_1}{n} \right)^{k_1} \cdots \left( \frac{j_m}{n}
        \right)^{k_m}.
\end{eqnarray*}
But
\begin{eqnarray*}
\lefteqn{\frac{1}{{n \choose {m - 1}}} \sum_{\jj} \left( \frac{j_1}{n}
          \right)^{k_1} \cdots \left( \frac{j_m}{n} \right)^{k_m}} \\
  &\to& (m - 1)! \int\!x_1^{k_1} \cdots x_{m - 1}^{k_{m - 1}}\,\left( 1 - x_1 -
          \cdots x_{m - 1} \right)^{k_m}\,dx_1 \cdots dx_{m - 1} \\
  & = & (m - 1)! \frac{\Gamma(k_1 + 1) \cdots \Gamma(k_m + 1)}{\Gamma(k + m)} =
          \frac{1}{{k \choose {k_1, \ldots, k_m}} {{k + m - 1} \choose 
{m - 1}}},
\end{eqnarray*}
where the above integral is over $(x_1, \ldots, x_{m - 1}) \in [0, 1]^{m - 1}$
with sum not exceeding unity.  Since the number of terms in $\sum^{**}_{\kk}$
is ${{k + m - 1} \choose {m - 1}} - m$, we therefore have
\begin{eqnarray*}
r_n(2k)
  &=& \frac{(2 k)!}{2^k k!} \sigma^{2 k} n^k \frac{{{k + m - 1} \choose {m - 1}}
        - m}{{{k + m - 1} \choose {m - 1}}} + o(n^k) \\
  &=& \frac{(2 k)!}{2^k k!} \sigma^{2 k} n^k \left[ 1 - \frac{m! \Gamma(k + 1)}
        {\Gamma(k + m)} \right] + o(n^k), \qquad k \geq 2.
\end{eqnarray*}
Similarly,
\[
r_n(2 k  - 1) = o\left(n^{k - (1 / 2)}\right), \qquad k \geq 2.
\]
Now part~(c) of the ATT implies~\eqref{even} and~\eqref{odd}.
\end{proof}

The following lemma lies at the heart of the proof of \refT{T:small1}.

\begin{lemma}
\label{L:rconds}
In the context of CLT~I, the sequence $(r_n)$ defined at~\eqref{rdef} satisfies
the conditions~\eqref{5.1}.
\end{lemma}

\begin{proof}
Clearly
\begin{equation}
\label{rdeft}
r_n = \frac{1}{{n \choose {m - 1}}} \sum_{\jj} \big[ t_n + \mut_{j_1} +
\cdots + \mut_{j_m} - \mut_n \big]^2, \qquad n \geq m - 1,
\end{equation}
with 
\begin{equation}
\label{mutdef}
\mut_n := \mu_n - \mu (n + 1),
\end{equation}
which is $o\left(\sqrt{n}\right)$ by~\eqref{asymean}.
Recall the inequality
\begin{equation}
\label{squares}
\left[ \sum_{i = 1}^k \xi_i \right]^2 \leq k \sum_{i = 1}^k \xi^2_i
\end{equation}
for real numbers $\xi_1, \ldots, \xi_k$.  Applying this to~\eqref{rdeft},
\begin{equation}
\label{rineq}
\frac{r_n}{m + 2} \leq t^2_n + \mut^2_n + \frac{m}{{n \choose {m - 1}}} \sum_{j
= 0}^{n - (m -  1)} {{n - 1 - j} \choose {m - 2}} \mut^2_j,
\end{equation}
from which~\eqref{5.1}(a) for $(r_n)$ is evident.

To establish the summability of $r_n / n^2$ we need only establish that of
$\mut^2_n / n^2$.  Indeed we can then use~\eqref{rineq} again, together
with~\eqref{bsimpler} and the estimate
\begin{eqnarray*}
\lefteqn{\sum_{n = m - 1}^{\infty} n^{-2} \frac{m}{{n \choose {m - 1}}} \sum_{j
           = 0}^{n - (m -  1)} {{n - 1 - j} \choose {m - 2}} \mut^2_j} \\
  &  = & m (m - 1) \sum_{j = 0}^{\infty} \mut^2_j \sum_{n = j + m - 1}^{\infty}
           \frac{(n - 1 - j)_{m - 2}}{n^2 (n)_{m - 1}} \\
  &\leq& m (m - 1) \sum_{j = 0}^{\infty} \mut^2_j \sum_{n = j + m - 1}^{\infty}
           n^{-3} \\
  &  = & O\left(\sum \frac{\mut^2_j}{j^2}\right) < \infty.
\end{eqnarray*}

To establish the summability of $\mut^2_n / n^2$, we recall
from~\eqref{muntilde} and~\eqref{Stirling} that
\begin{eqnarray}
\label{mu4}
\mut_n
  &=& O(n^{\ga - 1}) + t_n - \frac{1}{H_m - 1} (n + 1) \sum_{k = n}^{\infty}
        \frac{\hat{t}_k}{(k + 1) (k + 2)} \\
  & & \hspace{.1in} + \sum_{j = 2}^{m - 1} O\left(n^{\ga_j - 1} \sum_{k =
        0}^{n - 1} \frac{|\hat{t}_k|}{(k + 1)^{\ga_j}}\right), \nonumber
\end{eqnarray}
writing $\ga_j := \Re \gl_j$ (with $\ga = \ga_2 < 3 / 2$, since $m \leq 26$).
Using~\eqref{squares}, we need only establish the summability of $n^{-2}$ times
the square of each of the four terms on the right in~\eqref{mu4}.  The first of
these verifications is trivial, and the second was carried out
at~\eqref{bsimpler}.  For the third we apply the Cauchy--Schwarz inequality
[compare~\eqref{bcs}] to give
\begin{eqnarray*}
\left[ \sum_{k = n}^{\infty} \frac{\hat{t}_k}{(k + 1) (k + 2)} \right]^2
  &\leq& \frac{1}{n} \left[ \sum_{k = n}^{\infty} \sqrt{n} \frac{\sqrt{k}}{(k +
           1) (k + 2)} \frac{|t_k|}{\sqrt{k}} \right]^2 \\
  &  = & O\left(\frac{1}{n} \sum_{k = n}^{\infty} \sqrt{n} k^{- 3 / 2}
           \frac{t^2_k}{k}\right) = O\left(n^{- 1 / 2} \sum_{k = 
n}^{\infty} k^{-
           5 / 2} t^2_k\right),
\end{eqnarray*}
whence
\begin{eqnarray*}
\sum_n \left[ \sum_{k = n}^{\infty} \frac{\hat{k}_k}{(k + 1) (k + 2)} \right]^2
  &=& O\left(\sum_n n^{- 1 / 2} \sum_{k = n}^{\infty} k^{- 5 / 2} 
t^2_k\right) \\
  &=& O\left(\sum_k k^{- 5 / 2} t^2_k k^{1 / 2}\right) = O\left(\sum_k
        \frac{t^2_k} {k^2}\right) < \infty
\end{eqnarray*}
by~\eqref{bsimpler} again.

We pause here to note that when $m = 2$ the proof is finished here, and that up
to now we have used only~\eqref{bsimpler}, not the stronger
assumption~\eqref{small1}(b).

For our fourth and final verification, it suffices [again by
invoking~\eqref{squares}] to establish the summability of
\begin{equation}
\label{4th}
n^{2 \rho - 4} \left[ \sum_{k = 1}^{n - 1} \frac{|t_k|}{k^{\rho}} \right]^2
\end{equation}
for any real $\rho < 3 / 2$.  To do this, we break the sum into $\sum_{k <
n^{\gd}}$ and $\sum_{n^{\gd} \leq k < n}$ and once again
invoke~\eqref{squares}.  In the range $\sum_{k < n^{\gd}}$ we simply
use $t_k = O\left(\sqrt{k}\right)$ and note
\[
n^{2 \rho - 4} \left[ \sum_{k < n^{\gd}} O\left(k^{(1 / 2) - \rho}\right)
\right]^2 = O\left(n^{2 \rho - 4} \left(n^{\gd}\right)^{3 - 2 \rho}\right) =
O\left(n^{\tau}\right)
\]
with $\tau < -1$.  In the range $\sum_{n^{\gd} \leq k < n}$ we use
Cauchy--Schwarz again:
\begin{eqnarray*}
n^{2 \rho - 4} \left[\sum_{n^{\gd} \leq k < n} \frac{|t_k|}{k^{\rho}} \right]^2
  &=& n^{2 \rho - 4} n^{3 - 2 \rho} \left[ \sum_{n^{\gd} \leq k < n}
  \frac{k^{(1 
        / 2) - \rho}}{n^{(3 / 2) - \rho}}\,\,\frac{|t_k|}{k^{1 / 2}} 
\right]^2 \\
  &=& O\left(n^{-1} \sum_{n^{\gd} \leq k < n} \frac{k^{(1 / 2) - \rho}}{n^{(3 /
        2) - \rho}}\,\,\frac{t^2_k}{k}\right) = O\left( n^{-1} \max_{n^{\gd}
        \leq k < n} \frac{t^2_k}{k}\right),
\end{eqnarray*}
which is summable by assumption~\eqref{small1}(b).
\end{proof}

Now we proceed to our ``borderline small'' case.

\begin{proof}[Proof of CLT~II (\refT{T:small2})]
Again we use the method of moments together with asymptotic transfer results.
Given the similarity to the proof of CLT~I, we will be brief here.

Equation~\eqref{asymean2} simply repeats the asymptotic transfer
result~\eqref{slowremainder}.  As before, $(\sigma^2_n)$ satisfies the
recurrence~\eqref{eq:2.7} with $(b_n)$ replaced by $(r_n)$
of~\eqref{rdeft}--\eqref{mutdef}, where again $\mu := K_1 /  (H_m - 1)$ and
$r_j := 0$ for $0 \leq j \leq m - 2$.  Here the proofs diverge somewhat.  The
analogue of \refL{L:rconds} is \refL{L:rconds2} below.  Then the asymptotic
variance assertions of CLT~II follow immediately from \refT{T:more}(c).

If $\sum^{\infty} \frac{L^2(k)}{k} < \infty$, then from~\eqref{bingham} applied
to $L^2$ it follows that~\eqref{10.1} 
satisfies~\eqref{even}--\eqref{odd} for $k
= 1$.  Then higher moments are treated exactly as in the proof of CLT~I to
complete the proof of CLT~II.

If $\sum^{\infty} \frac{L^2(k)}{k} = \infty$, then one uses~\eqref{bingham},
\refT{T:more}(d), and induction to show that the moments~\eqref{10.1} satisfy
\begin{equation}
\label{even2}
\mut_n(2 k) \sim \frac{(2 k)!}{2^k k!} s^{2 k}(n), \qquad k \geq 1
\end{equation}
and
\begin{equation}
\label{odd2}
\mut_n(2 k - 1) = o\left(s^{2 k - 1}(n)\right), \qquad k \geq 1
\end{equation}
and thereby complete the proof of CLT~II.  We omit the details.
\end{proof}

The following cousin to \refL{L:rconds} was used in the proof of
\refT{T:small2}.

\begin{lemma}
\label{L:rconds2}
In the context of CLT~II, the sequence $(r_n)$ defined for $n \geq m - 1$ by
\[
r_n := \frac{1}{{n \choose {m - 1}}} \sum_{\jj} \big[ t_n + \mut_{j_1} + \cdots
+ \mut_{j_m} - \mut_n \big]^2
\]
satisfies
\[
r_n \sim (H_m - 1) \sigma^2 n L^2(n).
\]
\end{lemma}

\begin{proof}
By~\eqref{asymean2}, with $\gth := \rising{(3 / 2)}{m - 1} / 
(m! - \rising{(3 / 2)}{m - 1})$, we have
\begin{eqnarray*}
r_n
  &\sim& \frac{1}{{n \choose {m - 1}}} \sum_{\jj} \left[ n^{1 / 2} L(n) - \gth
           j_1^{1 / 2} L(j_1) - \cdots - \gth j_m^{1 / 2} L(j_m) + 
\gth n^{1 / 2}
           L(n) \right]^2 \\
  &\sim& n L^2(n)\,(m - 1)! \int\!\left[ (1 + \gth) - \gth \sum_{i = 1}^{m - 1}
           x_i^{1 / 2} - \theta \left( 1 - \sum_{i = 1}^{m - 1} x_i \right)^{1 /
           2} \right]^2\,dx_1 \cdots dx_{m - 1},
\end{eqnarray*}
where the integral, call it~$J$, is over $(x_1, \ldots, x_{m - 1}) 
\in [0, 1]^{m
- 1}$ with sum not exceeding unity.  To complete the proof we need only show
$J = (H_m - 1) \sigma^2 / (m - 1)!$, with $\sigma^2$ defined
at~\eqref{sigmadef2}.

Indeed,
\begin{eqnarray*}
J
  &=& \int\!\left[ (1 + \gth)^2 + \gth^2 \sum_{i = 1}^{m - 1} x_i + \theta^2
        \left( 1 - \sum_{i = 1}^{m - 1} x_i \right) \right. \\
  & & \qquad {} - 2 \gth (1 + \gth) \sum_{i = 1}^{m - 1} x_i^{1 / 2} - 2 \gth (1
        + \gth) \left( 1 - \sum_{i = 1}^{m - 1} x_i \right)^{1 / 2} \\
  & & \qquad \left. {} + \gth^2 \sum_{i, j:\,i \neq j} x_i^{1 / 2} x_j^{1 / 2} +
        2 \gth^2 \left( 1 - \sum_{i = 1}^{m - 1} x_i \right)^{1 / 2} \sum_{i =
        1}^{m - 1} x_i^{1 / 2} \right]\,dx_1 \cdots dx_{m - 1} \\
  &=& \left[ (1 + \gth)^2 + \gth^2 \right]\,\frac{1}{(m - 1)!} - 2 \gth (1 +
        \gth) m \frac{\Gamma(3 / 2)}{\Gamma(m + (1 / 2))} + \gth^2 m (m - 1)
        \frac{[\Gamma(3 / 2)]^2}{\Gamma(m + 1)} \\
  &=& \left[ (1 + \gth)^2 + \gth^2 \right]\,\frac{1}{(m - 1)!} - 2 \gth (1 +
        \gth) m\,\frac{1}{\rising{(3 / 2)}{m - 1}} + \gth^2 \frac{\pi /
        4}{(m - 2)!}.
\end{eqnarray*}
Plugging in the value of~$\gth$ and simplifying, we obtain $J = (H_m - 1)
\sigma^2 / (m - 1)!$, as desired.

\end{proof}

\subsection{Periodicity for $m \geq 27$}
\label{S:phase}

If $t_n = o\left(\sqrt{n}\right)$ as in CLT~I but $m \geq 27$, then the
remainder term $\mut_n := \mu_n - \mu (n + 1)$ for the mean---which
by~\eqref{asymean} was $o\left(\sqrt{n}\right)$ when $m \leq 26$---now
satisfies, by the ETT and~\eqref{exactlinear2} [compare~\eqref{muntilde}]
\begin{eqnarray}
  \label{eq:peridocity}
\mut_n
  &=& c_2 \frac{n^{\gl_2 - 1}}{\Gamma(\gl_2)} + c_3 \frac{n^{\gl_3 -
        1}}{\Gamma(\gl_3)} \notag\\
  & & \hspace{.1in} + m! \sum_{j = 2}^{m - 1} \frac{1}{\psi'(\gl_j)} 
[z^n] \left(
        (1 - z)^{- \gl_j} \int_{0}^z\!\hat{T}(\zeta) (1 - \zeta)^{\gl_j -
        1}\,d\zeta \right) \notag\\
  & & \hspace{.1in} + o\left(\sqrt{n}\right) + O\left(n^{\Re \gl_4 - 1}\right).
\end{eqnarray}

\emph{Typically} this will lead to the negative result that $(\mut_n)$ [and
hence also $(r_n)$ and $(\sigma^2_n)$] suffers from periodicity and that there
is no natural distributional limit for normalized~$X_n$.  Examples are
the space requirement functional 
studied by Chern and Hwang~\cite{MR1871558} (and others before
them~\cite{MR90a:68012})
and the shape functional~\cite{MR97f:68021}.  For recent developments
in the case~$m > 26$ for the space requirement
see~\cite{Chauvin-Pouyanne-2003}.

But it is perhaps difficult
to establish a \emph{general} negative result, due to
cancellations.  For example, suppose $T(z)$ equals $(1 - z)^{-1}$, so that
$t_n \equiv 1$, as studied by Chern and Hwang~\cite{MR1871558}, except
perhaps that the initial
values $t_0, \ldots, t_{m - 2}$ are changed.  Then
\[
T^{(m - 1)}(z) \equiv (m - 1)! (1 - z)^{- m},
\]
whence
\begin{eqnarray*}
A(z)
  &=& \sum_{n=1}^\infty\mu_n z^n = \sum_{j = 1}^{m - 1} c_j (1 - z)^{-
  \gl_j} + (m -  1)! \sum_{j =
        1}^{m - 1} \frac{(1 - z)^{- \gl_j}}{\psi'(\gl_j)} \int_{0}^z\!(1 
        - \zeta)^{\gl_j - 2}\,d\zeta \\
  &=& \sum_{j = 1}^{m - 1} c_j (1 - z)^{- \gl_j} - \frac{1}{m - 1} (1 - z)^{-1}
        - (m - 1)! \sum_{j = 1}^{m - 1} \frac{(1 - z)^{- \gl_j}}{(1 -
        \gl_j)\psi'(\gl_j)}.
\end{eqnarray*}
Now it is possible to choose $t_0, \ldots, t_{m - 2}$ so that
\begin{equation}
\label{ctrex}
c_j = \frac{(m - 1)!}{(1 - \gl_j) \psi'(\gl_j)}, \qquad j = 1, \ldots, m - 1.
\end{equation}
In that case $A(z) = - \frac{1}{m - 1} (1 - z)^{-c}$, whence $\mut_n \equiv
\mu_n \equiv - 1/ (m - 1)$  [and we see that the
chosen values of $t_0, \ldots, t_{m - 2}$ are all $ - 1 / (m - 1)$], and we get
linear variance and asymptotic normality, just as in CLT~I, for \emph{every} $2
\leq m < \infty$.

One might object that the above example is artificial, in that the toll
sequence changes sign.  But the same calculation show that if the toll sequence
is chosen as above ($t_n \equiv 1$) but with inital values%
\[
t_j := K (j + 1) - \sfrac{1}{m - 1}, \qquad 0 \leq j \leq m - 2,
\]
then still, for every $m \geq 2$, the sequence $(\mut_n)$ is constant, the
variance is linear, and we have asymptotic normality.  Further, $(t_n)$ is
nonnegative provided $K \geq 1 / (m - 1)$.  [We remark in passing that the
choice $K = 1 / (m - 1)$ leads to the degenerate case of \refR{R:degenerate}.]
The sequence $(t_n)$ is also nondecreasing (as in most real examples) provided
$K \leq m / (m - 1)^2$.



\section{Moderate and large toll functions}
\label{sec:moderate-large-toll}


In order to describe the limiting distribution of $X_n$ for moderate
and large tolls, we will introduce
a family of random variables $Y \equiv Y(\beta)$ defined for
$\beta > 1/2$, $\beta \ne 1$.  Anticipating Lemma~\ref{lemma:6.1.1},
we need to consider the distributional equation
\begin{equation}
  \label{eq:1}
  Y \stackrel{\mathcal{L}}{=} \sum_{j=1}^{m} S_j^\beta
  Y_{j} + 1
\end{equation}
Here $(Y_j)_{j=1}^m$ are independent copies of $Y$ and $(S_1,\ldots,S_m)$
is uniformly distributed on the $(m-1)$-simplex,
independent of
$(Y_j)_{j=1}^m$.  Recall that the $(m-1)$-simplex is the set
\[
\{(s_1,\ldots,s_{m})\!: s_j \geq 0 \text{ for } 1 \leq j \leq m
\text{ and } \mathbf{s_+} = 1 \}, \]
where $\mathbf{s_+}$ denotes $\sum_{j=1}^{m} s_j$.

Let
$U_{(1)},\ldots,U_{(m-1)}$ be the order statistics of a sample of size $m-1$
from the uniform distribution on $(0,1)$. They have joint density
\begin{equation*}
f_{U_{(1)},\ldots,U_{(m-1)}}(x_1,\ldots,x_{m-1}) \equiv (m-1)!
\indicator{0 < x_1 <  \cdots < x_{m-1} < 1}
\end{equation*}
with respect to Lebesgue measure on $\mathbb{R}^{m-1}$, 
where $\indicator{A}$ is the indicator of $A$.

{\sloppy
By a change of variables,  we find that the
joint distribution of the spacings $S_1,\ldots,S_{m}$, defined, with
$U_{(0)} := 0$ and $U_{(m)} := 1,$  by
\begin{equation*}
  S_{i} := U_{(i)} - U_{(i-1)}, \; i=1,\ldots,m,
\end{equation*}
is uniform over the $(m-1)$-simplex:
\begin{equation*}
  f_{S_1,\ldots,S_{m-1}}(s_1,\ldots,s_{m-1}) \equiv (m-1)! \indicator{s_j > 0,
    j=1,\ldots,m-1; \; \mathbf{s_{+}} < 1}.
\end{equation*}%
}
When $r_j > -1$ for $1 \leq j \leq m$, observe that
\begin{align}
  \EE\,\left[\prod_{j=1}^m S_j^{r_j}\right] & = (m-1)!
  \int_{\substack{s_1,\ldots,s_{m-1} > 0\\ \mathbf{s_{+}}<1}}
  s_1^{r_1} \cdots s_{m-1}^{r_{m-1}}(1 - \mathbf{s_{+}})^{r_m} \,
  ds_{m-1}\cdots ds_1\notag\\
  \label{eq:2}
  & =: (m-1)! B(r_1+1,\ldots,r_m + 1) \quad \text{(defining $B$ as the
  integral)} \\
  & = (m-1)! \frac{\prod_{j=1}^{m}
  \Gamma(r_j + 1)}{\Gamma(r_1+\cdots+r_m + m)}. \notag 
\end{align}

\begin{lemma}
  \label{lemma:6.1.1}
  Fix $\beta > 1/2$ with $\beta \ne 1$.  Then there
  exists a unique distribution $\mathcal{L}(Y) \equiv
  \mathcal{L}(Y(\beta))$ with finite second
  moment satisfying the distributional identity~\eqref{eq:1}.
\end{lemma}
\begin{proof}
  We first observe the that mean of any such distribution is
  determined by~\eqref{eq:1}.  Indeed, by taking expectations
  in~\eqref{eq:1} and using~\eqref{eq:2}, we get
  \begin{equation*}
    \mu := \EE\,Y = \left( 1 - \frac{m!
        \Gamma(\beta+1)}{\Gamma(\beta+m)} \right)^{-1}
  \end{equation*}
  since $\beta \ne 1$.  Thus we can equivalently
  consider the distributional identity
  \begin{equation*}
    W \stackrel{\mathcal{L}}{=} \sum_{j=1}^m S_j^\beta W_j + H,
  \end{equation*}
  where
  \begin{equation*}
    H := 1 - \mu + \mu \sum_{j=1}^m S_j^\beta.
  \end{equation*} Here $W$ is restricted to have mean~0 and finite
  second moment, $(W_j)_{j=1}^m$ are independent copies of $W$, and
  $(S_1,\ldots,S_m)$ is uniformly distributed on the $(m-1)$-simplex,
  independent of $(W_j)_{j=1}^m$.

  We now employ a standard contraction-method
  argument~\cite{MR2003b:68217,MR2003c:68096}.  Let $d_2$
  denote the metric on $\mathcal{M}_2(0)$, the space of probability
  distributions with mean~0 and finite variance, defined by
  \begin{equation*}
    d_2( G_1, G_2) := \min \lVert X_2 - X_1 \rVert_2,
  \end{equation*}
  taking the minimum over all pairs of random variables $X_1$ and
  $X_2$ defined on a common probability space with $\mathcal{L}(X_1)
  = G_1$ and $\mathcal{L}(X_2) = G_2$.  Here $\lVert \cdot \rVert_2$
  denotes  $L_2$-norm.

  Let $T$ be the map
  \begin{equation*}
    T: \mathcal{M}_2(0) \to \mathcal{M}_2(0), \quad G \mapsto
    \mathcal{L} \left( \sum_{j=1}^{m} S_j^\beta X_j + H \right),
  \end{equation*}
  where $(X_j)_{j=1}^m$ are independent with $\mathcal{L}(X_j) = G$,
  $j=1,\ldots,m$, and $(S_1,\ldots,S_m)$ is uniformly distributed on
  the $(m-1)$-simplex, independent of $(X_j)_{j=1}^m$.  We show that
  $T$ is a contraction on
  $\mathcal{M}_2(0)$; more precisely, that there exists a $\rho < 1$
  such that
  \begin{equation*}
    d_2(T(\mathcal{L}(A)),T(\mathcal{L}(B))) \leq \rho
    d_2(\mathcal{L}(A),\mathcal{L}(B))
  \end{equation*}
  for all pairs $\mathcal{L}(A)$ and $\mathcal{L}(B)$ in
  $\mathcal{M}_2(0)$.  To bound
  $d_2(T(\mathcal{L}(A)),T(\mathcal{L}(B)))$, we couple $T(\mathcal{L}(A))$ and
  $T(\mathcal{L}(B))$ by taking $m$ independent copies  $(A_j, B_j)$
  of the $d_2$-optimally coupled $(A,B)$,  an independent
  $(S_1,\ldots,S_m)$, and defining
  \begin{equation*}
    A' := \sum_{j=1}^m S_j^\beta A_j + H \sim T(\mathcal{L}(A)),\quad
    B' := \sum_{j=1}^m S_j^\beta B_j + H \sim T(\mathcal{L}(B)).
  \end{equation*}
  Now, defining $\mathbf{S} := (S_1,\ldots,S_m)$ and  using the law of
  total variance, 
  \begin{align*}
    &d_2( \mathcal{L}(T(A)), \mathcal{L}(T(B)) )^2\\
    \quad &\leq \lVert B' - A'
    \rVert_2^2
    = \left \lVert \sum_{j=1}^m S_j^\beta (B_j - A_j)
    \right\rVert_2^2
    = \Var\, \left[ \sum_{j=1}^m S_j^\beta (B_j - A_j) \right]\\
    &= \EE\, \Var\, \left[ \left.\sum_{j=1}^m S_j^\beta (B_j - A_j) \right|
      \mathbf{S} \right]
    + \Var\,
    \EE\, \left[ \left. \sum_{j=1}^m S_j^\beta (B_j - A_j) \right|
    \mathbf{S} \right]\\
    &= \sum_{j=1}^m (\EE S_j^{2\beta}) \Var\, [ B_j - A_j ] =
    d_2(\mathcal{L}(A),\mathcal{L}(B))^2  \sum_{j=1}^m \EE S_j^{2\beta}
    = m! \frac{\Gamma(2\beta+1)}{\Gamma(2\beta+m)}
    d_2(\mathcal{L}(A),\mathcal{L}(B))^2.
  \end{align*}
  We need only verify that
  \begin{equation*}
    \rho^2 := m! \frac{\Gamma(2\beta+1)}{\Gamma(2\beta+m)} =
    \frac{m!}{(2\beta+m-1) \cdots (2\beta+1)} < 1,
  \end{equation*}
  which is true when $\beta > 1/2$.  The existence and uniqueness of
  $\mathcal{L}(Y)$ now follows from the Banach fixed point
  theorem~\cite[Theorem~2]{MR2003c:68096}.
\end{proof}

\subsection{Moderate toll functions}
\label{sec:moder-toll-funct}
In the case of moderate toll functions, convergence in distribution
and convergence of all moments can be stated as
\begin{theorem}
  \label{thm:F2}
  If the toll sequence~$(t_n)$ satisfies
\[
t_n \sim n^{\beta}L(n) \text{ with $1/2 < \beta < 1$},
\]
where $L$ is a slowly varying function
and $\alpha < 1 + \beta$, then the mean of the corresponding additive
functional $X_n$ on $m$-ary search trees with the random permutation
model satisfies
\begin{equation}
  \label{eq:F2.1}
  \mu_n = \mu n -
  \frac{\rising{(1+\beta)}{m-1}}{m!-\rising{(1+\beta)}{m-1}}
  n^{\beta}L(n) + o(n^{\beta} L(n)),  \qquad \mu := \frac{K_1}{H_m-1},
\end{equation}
with $K_1$ defined at~\eqref{5.1a}.  Moreover,
\[
\frac{X_n-\mu{}n}{n^{\beta} L(n)} \stackrel{\mathcal{L}}{\to}
Y_{\beta},
\]
with convergence of all moments.
\end{theorem}
\begin{remark}
  \label{remark:alpha-monotone}
It is well known that $\alpha < 3/2$ for $m \leq 26$.  In
Theorem~\ref{theorem:alpha-monotone} we show that 
$\alpha$ increases with $m$.  Thus for a fixed $\beta \in (1/2,1)$,
the condition $\alpha < 1 + \beta$  is
equivalent to $m \leq m_0$ for some $m_0 \geq 26$.
\end{remark}
\begin{proof}[Proof of Theorem~\ref{thm:F2}]
  We use the notation introduced in the proof of
  Theorem~\ref{T:small1} in Section~\ref{S:CLTpfs}.
  Equation~(\ref{eq:F2.1}) is simply a restatement of the asymptotic
  transfer result~\eqref{slowremainder}.

We show that the moments $\tilde\mu_n(k)$ satisfy
\begin{equation}
  \label{eq:F2.5}
  \tilde\mu_n(k) = g_k n^{k\beta}L^k(n) + o(n^{k\beta} L^k(n)) \text{
  as } n   \to \infty.
\end{equation}
The claim holds for $k=1$ with
\begin{equation}
  \label{eq:8}
g_1 := -\frac{\rising{(1+\beta)}{m-1}}{m!-\rising{(1+\beta)}{m-1}} =
\left( 1 - \frac{m!\Gamma(\beta+1)}{\Gamma(\beta+m)} \right)^{-1}.
\end{equation}
Using~\eqref{eq:11}, by induction we get, for $k \geq 2$,
\begin{align*}
  r_n(k) &= o(n^{k\beta}L^k(n)) \\
  & \quad {} +  \sideset{}{_{}^*}\sum_{\mathbf{k}}
  \binom{k}{k_1,\ldots,k_m,k_{m+1}}
  (n^{\beta}L(n))^{k_{m+1}} \\
  &\qquad {} \times \frac{1}{\binom{n}{m-1}}
  \sum_{\mathbf{j}} g_{k_1}(j_1^\beta L(j_1))^{k_1} \cdots
  g_{k_m}(j_m^\beta L(j_m))^{k_m}\\
  &= o(n^{k\beta} L^k(n))\\
  & \quad {}+ \sideset{}{_{}^*}\sum_{\mathbf{k}}
  \binom{k}{k_1,\ldots,k_m,k_{m+1}}
  (n^\beta L(n))^{k} g_{k_1} \cdots g_{k_m} \\
  & \qquad {} \times \frac{1}{\binom{n}{m-1}}
  \sum_{\mathbf{j}} \left(\frac{j_1}{n}\right)^{k_1\beta}
  \cdots \left(\frac{j_m}{n}\right)^{k_m\beta} \frac{L^{k_1}(j_1)
  \cdots L^{k_m}(j_m)}{L^{k_1 + \dots + k_m}(n)}.
\end{align*}
But [recall the definition of $B$ at~\eqref{eq:2}]
\[
\frac{1}{\binom{n}{m-1}} \sum_{\mathbf{j}}
\left(\frac{j_1}{n}\right)^{k_1\beta}  \cdots
\left(\frac{j_m}{n}\right)^{k_m\beta} \frac{L^{k_1}(j_1)
  \cdots L^{k_m}(j_m)}{L^{k_1 + \dots + k_m}(n)}
\to (m-1)!
B(k_1\beta+1,\ldots k_m\beta+1)
\]
so that
\begin{multline*}
    r_n(k) = o(n^{k\beta}L^k(n)) \\
    +  n^{k\beta} L^k(n)
    (m-1)!
    \sideset{}{_{}^*}\sum_{\mathbf{k}} 
    \binom{k}{k_1,\ldots,k_m,k_{m+1}} g_{k_1}\cdots g_{k_m}
    B(k_1\beta+1,\ldots,k_m\beta+1).
\end{multline*}
Using Theorem~\ref{T:more}, with $v=k\beta > 1$, we get
\begin{multline*}
  \tilde\mu_n(k) = o(n^{k\beta} L^k(n))\\
  +  n^{k\beta} L^k(n) \frac{(m-1)!}{1-\frac{m!
    \Gamma(k\beta+1)}{\Gamma(k\beta+m)}}
    \sideset{}{_{}^*}\sum_{\mathbf{k}}
  \binom{k}{k_1,\ldots,k_m,k_{m+1}} g_{k_1}
    \cdots g_{k_m} B(k_1\beta+1,\ldots,k_m\beta+1).
\end{multline*}
Thus, defining $g_k$ recursively as
\begin{equation}
  \label{eq:9}
g_k =
\frac{(m-1)!}{1-\frac{m!\Gamma(k\beta+1)}{\Gamma(k\beta+m)}}
\sideset{}{_{}^*}\sum_{\mathbf{k}} \binom{k}{k_1,\ldots,k_{m+1}}
g_{k_1}\cdots g_{k_m} B(k_1\beta+1,\ldots,k_m\beta+1),
\end{equation}
with $g_0 = 1$,
we see that~\eqref{eq:F2.5} holds for all $k \geq 0$.

By Lemma~\ref{lem:F7.1} (to follow) and the method of moments
(cf., e.g.,~\cite[Sections~4.4 and~4.5]{MR49:11579}), the $g_k$'s are the
moments of a uniquely determined distribution, say
$\mathcal{L}(\widehat{Y})$, and
\begin{equation*}
  \frac{X_n - \mu n}{n^\beta L(n)} \stackrel{\mathcal{L}}{\to} \widehat{Y}
\end{equation*}
with convergence of all moments.  It remains to show that $\widehat{Y}
\stackrel{\mathcal{L}}{=} Y(\beta)$.

Define
\begin{equation}
  \label{eq:7}
  \widetilde{Y} := \sum_{j=1}^m S_j^\beta \widehat{Y}_j + 1,
\end{equation}
where $(\widehat{Y}_j)_{j=1}^m$ are independent copies of $\widehat{Y}$ and
$(S_1,\ldots,S_m)$ is uniformly distributed on the $(m-1)$-simplex,
independent of $(\widehat{Y}_j)_{j=1}^m$.  We will show that $\widehat{Y}
\stackrel{\mathcal{L}}{=} \widetilde{Y}$, and then, by~\eqref{eq:7},
$\mathcal{L}(\widehat{Y})$ satisfies the distributional
identity~\eqref{eq:1} and has finite second moment.  By
Lemma~\ref{lemma:6.1.1}, $\widehat{Y} \stackrel{\mathcal{L}}{=}
Y(\beta)$, as desired.

To show $\widehat{Y} \stackrel{\mathcal{L}}{=} \widetilde{Y}$, it
suffices to show that
$\widehat{Y}$ and $\widetilde{Y}$ have the same moments.  Letting
$\sum_{\mathbf{k}}$ denote (as before) the sum over $(m+1)$-tuples
$(k_1,\ldots,k_{m+1})$ of nonnegative integers summing to $k$, and
using~\eqref{eq:7},
\eqref{eq:2}, and~\eqref{eq:9},
\begin{align*}
  \EE\,\widetilde{Y}^k &= \sum_{\mathbf{k}} \binom{k}{k_1,\ldots,k_{m+1}}
  \EE\,\left[\prod_{j=1}^m (S_j^{\beta} \widehat{Y}_j)^{k_j} \right] \\
  &= \sum_{\mathbf{k}} \binom{k}{k_1,\ldots,k_{m+1}} (m-1)! B(k_1\beta
  + 1, \ldots k_m\beta + 1) g_{k_1} \cdots g_{k_m}\\
  &= (m-1)! \sideset{}{_{}^*}\sum_{\mathbf{k}}
  \binom{k}{k_1,\ldots,k_{m+1}}  B(k_1\beta 
  + 1, \ldots k_m\beta + 1) g_{k_1} \cdots g_{k_m}\\
  & \qquad + m!  B(k\beta+1,1,\ldots,1) g_k\\
  &= \left[ 1 - \frac{m! \Gamma(k\beta+1)}{\Gamma(k\beta+m)} \right]
  g_k + \frac{m! \Gamma(k\beta+1)}{\Gamma(k\beta+m)} g_k = g_k = \EE\,
  \widehat{Y}^k,
\end{align*}
where (as before) $\sum_{\mathbf{k}}^{\mathbf{*}}$ denotes the same
sum as $\sum_{\mathbf{k}}$ with the additional restriction that $k_i <
k$ for $i=1,\ldots,m$.
\end{proof}

\begin{lemma}
  \label{lem:F7.1}
  The moments $(g_k)$ uniquely determine the
  distribution $\mathcal{L}(Y)$.
\end{lemma}
\begin{proof}
  Define $\gamma_k := g_k/k!$.  It suffices to show (by Carleman's
  condition) that there
  exists an $M$ such that $\gamma_k \leq M^k$ for all $k$ sufficiently
  large.  We proceed by induction.  Indeed, by~\eqref{eq:9} we know
  \begin{align*}
    \gamma_k &= \frac{(m-1)!}{1 -
    \frac{m!\Gamma(k\beta+1)}{\Gamma(k\beta+m)}}
    \sideset{}{_{}^*}\sum_{\mathbf{k}}
    \frac{1}{k_{m+1}!} \left( \prod_{j=1}^m \gamma_{k_j} \right)
    B(k_1\beta+1, \ldots, k_m\beta+1)\\
    & \leq M^k \frac{(m-1)!}{1 -
    \frac{m!\Gamma(k\beta+1)}{\Gamma(k\beta+m)}}
    \sideset{}{_{}^*}\sum_{\mathbf{k}}
    \frac{M^{-k_{m+1}}}{k_{m+1}!} B(k_1\beta+1,\ldots,k_m\beta+1)
  \end{align*}
  by the induction hypothesis.  So it is certainly sufficient to show
  that
  \begin{align}
    \label{eq:3}
    &\sideset{}{_{}^*}\sum_{\mathbf{k}}
    \frac{M^{-k_{m+1}}}{k_{m+1}!} B(k_1\beta+1,\ldots,k_m\beta+1)
    \notag\\ 
    &= \sum_{k_{m+1}=0}^k \frac{M^{-k_{m+1}}}{k_{m+1}!} \Gamma(
    (k-k_{m+1})\beta + m)^{-1} \sum_{\substack{0 \leq k_1,\ldots,k_m <
    k\\ k_1 + \cdots + k_m = k - k_{m+1}}} \prod_{j=1}^m \Gamma(
    k_j\beta + 1) \notag\\
    &\to    0 \quad \text{ as $k \to \infty$}.
  \end{align}
  
  For this, fix a value of $k_{m+1} \in \{0,1,2,\ldots\}$, and
  consider the sum
  \begin{equation}
    \label{eq:5}
     \sum_{\substack{0 \leq
    k_1,\ldots,k_m < k\\ k_1+\cdots+k_m=k-k_{m+1}}} \prod_{j=1}^m
    \Gamma(k_j\beta + 1).
  \end{equation}
  By log-convexity of $\Gamma$~\cite[6.4.1]{MR29:4914},
  taking $I$ to be $(0,\infty)$ and $g$ to be $\log \Gamma$ in
  Proposition~3.C.1 of~\cite{MR81b:00002}, the logarithm of the
  product in~\eqref{eq:5} is Schur-convex on $(0,\infty)^m$.  Thus applying
  Proposition~5.C.2 of~\cite{MR81b:00002} with $m=0$ there, the
  biggest terms in
  the sum correspond to $k_{j} = k-k_{m+1}$ for $j$ equal to some
  $j_0$ and $k_j=0$ otherwise; together, these $m$ terms contribute
  $m\Gamma((k-k_{m+1})\beta + 1)$ to the sum.  If $k > k_{m+1}$, there
  are other terms, the biggest of which corresponds to having one of
  the $k_j$'s be $k-k_{m+1}-1$, one be~1, and the rest be~0.  (This
  follows from Proposition~5.C.1 of~\cite{MR81b:00002} with $m=0$ and
  $M=k-k_{m+1}-1$ there.)  The  total number of terms in the
  sum~\eqref{eq:5} is at most   $\binom{m-1+(k-k_{m+1})}{m-1}$.  So
  the remaining  contribution  to~\eqref{eq:5} is at most
  \begin{equation*}
    {\binom{m-1+(k-k_{m+1})}{m-1}}
    \Gamma( (k-k_{m+1}-1)\beta + 1) \Gamma(\beta + 1).
  \end{equation*}
  We have found that the left side of~\eqref{eq:3} is bounded by
  \begin{equation*}
    \sum_{k_{m+1}=0}^{\infty} \frac{M^{-k_{m+1}}}{k_{m+1}!}
    \indicator{k_{m+1} \leq k} f(k-k_{m+1}),
  \end{equation*}
  where
  \begin{align}
    \label{eq:6}
    f(k) & \leq \frac{m}{(k\beta+1) \cdots (k\beta+m-1)} +
    \Gamma(\beta+1) \frac{\binom{m-1+k}{m-1}\Gamma((k-1)\beta +
      1)}{\Gamma(k\beta+m)} \\
    & \leq \frac{\binom{m-1+k}{m-1}\Gamma(k\beta+1)}{\Gamma(k\beta+m)}
    = \frac{1}{(m-1)!} \frac{(k+1) \cdots (k+(m-1))}{(k\beta+1) \cdots
    (k\beta + (m-1))},
  \end{align}
  which is a bounded function of $k$.  To apply the dominated
  convergence theorem, it suffices to show that the right side
  of~\eqref{eq:6} tends to 0 as $k \to \infty$, which follows from
  Stirling's approximation and the fact that $\beta > 0$.
\end{proof}

When $t_n$ satisfies the conditions in Theorem~\ref{thm:F2} but
$\alpha \geq 1+\beta$ then [compare~\eqref{eq:peridocity}],%
\begin{align*}
  \tilde\mu_n &= c_2\frac{n^{\lambda_2-1}}{\Gamma(\lambda_2)} + c_3
  \frac{n^{\lambda_3-1}}{\Gamma(\lambda_3)} \\
  & \qquad + m!\sum_{j=1}^{m-1} \frac{1}{\psi'(\lambda_j)} [z^n]
  \left((1-z)^{-\lambda_j} \int_{0}^z
  \hat{T}(\zeta)(1-\zeta)^{\lambda_j-1}\,d\zeta\right) \\
  & \qquad + o(n^{\beta}) + O(n^{\Re(\lambda_4-1)}),
\end{align*}
and typically this leads to periodicity.

\subsection{Large toll functions}
\label{sec:large-toll-functions}

If $t_n \sim n^\beta L(n)$, where $\beta > 1$ and $L$ is slowly
varying function, then we have convergence in
distribution for all values of $m$.  We state the result, omitting the
proof, as it is very similar to that of Theorem~\ref{thm:F2}.
\begin{theorem}
  \label{thm:F3}
  If the toll sequence~$(t_n)$ satisfies
  \[
  t_n \sim n^\beta L(n), \text{with $\beta > 1$},
  \]
  where $L$ is a slowly varying function,
  then
  \[
  \frac{X_n}{n^\beta L(n)} \stackrel{\mathcal{L}}{\to} Y_(\beta)
  \]
  with convergence of all moments, where $\mathcal{L}(Y(\beta))$ is the
  unique distribution satisfying~\eqref{eq:1}.
\end{theorem}

Presumably, the borderline case $t_n \sim n L(n)$ where $L$ is a slowly
varying function can also be handled using the techniques of this
paper, but we have not pursued this.  The specific choice $t_n \equiv
n-(m-1)$ for $n \geq m-1$ corresponds to the well-studied total path
length of a
random $m$-ary search tree.  The corresponding additive functional
measures the number of basic operations in $m$-ary \texttt{Quicksort}.
As is well known, the number of basic operations has
mean $\Theta(n \log n)$ and standard deviation $\Theta(n)$.
See~\cite{hennequin91:_analy} for details in this case
and~\cite[Corollary~5.2]{MR2000h:60034} for a characterization of the
limiting distribution of the path length.




\appendix
\section{Solution of an Euler differential equation}
\label{appendix:diffeq}

We now provide the proof of equation~\eqref{1.2e} in
Theorem~\ref{T:ett}, which states
the general solution of the differential equation~\eqref{eq:2.8} with
initial conditions \( a_j = b_j \), \( 0 \leq j \leq m-2 \).  This
linear differential equation can be written in the form
\begin{equation}
  \label{eq:A.1}
  \mathbf{L}g = h
\end{equation}
where the operator \( \mathbf{L} \) is defined as
\begin{equation}
  \label{eq:A.1.5}
  (\mathbf{L}g)(z) := g^{(m-1)}(z) - m! (1-z)^{-(m-1)} g(z).
\end{equation}
We seek the solution $g=A$ corresponding to input
$h = B^{(m-1)}$.

Equations of the form~\eqref{eq:A.1}--\eqref{eq:A.1.5} are members of a
class known as \emph{Euler differential equations}. In this appendix
we discuss a general method for solving Euler equations, restricting
attention, for the sake of definiteness and practicality,
to~\eqref{eq:A.1}--\eqref{eq:A.1.5}.  We assume that the reader is
familiar with the theory of linear differential equations with
constant coefficients (see,
e.g.,~\cite{boyce86:_elemen_differ_equat}).

For brevity we have omitted several routine proofs.  A fuller version of
Appendices~\ref{appendix:diffeq} and~\ref{appendix:indicial} may be
found in the technical report~\cite{FK-transfer-arXiv}.

\subsection*{The homogeneous solution}
\label{sec:homogeneous-solution}

The technique for solving \( \mathbf{L}g = 0 \) is quite easily
summarized: make the change of variable \( z = 1 - e^{-x} \), that is,
\( x = \ln{((1-z)^{-1}}) \).  For notational convenience we will
abbreviate $\ln{((1-z)^{-1}})$ as $L(z)$.  Then consider the function
\( \tilde{g} \) defined by
\begin{equation}
  \label{eq:A.2}
  \tilde{g}(x) := g(1-e^{-x}), \quad \text{ i.e., } \quad g(z) =
  \tilde{g}(L(z)).
\end{equation}
\begin{lemma}
  \label{lemma:A.1}
  The derivatives of \( g \) are related to those of \( \tilde{g} \)
  by
  \begin{equation}
    \label{eq:A.3}
    {g}^{(k)}(z) = (1-z)^{-k} \sum_{j=0}^k \stirlingone{k}{j}
    \tilde{g}^{(j)}(L(z)),
  \end{equation}
  where \( \stirlingone{k}{j} \) denotes a signless Stirling number
  of the first kind.
\end{lemma}
\begin{proof}
  The proof is a straightforward induction on~$k$ using standard
  identities for the Stirling numbers~\cite[\S~1.2.6]{knuth97}. 
\end{proof}
The left-hand side of~\eqref{eq:A.1} can hence be expressed as
\begin{equation*}
  (\mathbf{L}g)(z) = (1-z)^{-(m-1)} \left\{ \sum_{j=0}^{m-1}
  \stirlingone{m-1}{j} \tilde{g}^{(j)}(L(z)) - m!
  \tilde{g}(L(z)) \right\}
\end{equation*}
so that solving \( \mathbf{L}g = 0 \) is equivalent to solving \(
\mathbf{\widetilde{L}}\tilde{g} = 0 \), where
\begin{equation}
  \label{eq:A.4}
  \mathbf{\widetilde{L}}\tilde{g}(x) := \sum_{j=0}^{m-1} \stirlingone{m-1}{j}
  \tilde{g}^{(j)}(x) - m! \tilde{g}(x).
\end{equation}
But this is a linear differential equation with constant coefficients.
Its \emph{indicial polynomial}, or \emph{characteristic polynomial},
is
\begin{equation}
  \label{eq:A.5}
  \psi_m(\lambda) \equiv \psi(\lambda) := \sum_{j=0}^{m-1}
  \stirlingone{m-1}{j} \lambda^j - m! = \rising{\lambda}{m-1} - m!,
\end{equation}
the last equality following from~\cite[1.2.6-(44)]{knuth97}.  For more
on this polynomial see Appendix~\ref{sec:prop-indic-polyn}.  From the
discussion in~\cite{MR93f:68045} it follows that there are \( m - 1\)
distinct (and nonzero) roots of \( \psi \), call them \(
\lambda_1,\ldots, \lambda_{m-1} \) arranged in nonincreasing order of
real parts.  Thus the functions \( \exp{( \lambda_j x )} \) are \( m
-1 \) linearly independent solutions of~\eqref{eq:A.4} and hence the
functions \( (1-z)^{-\lambda_j} \) form a basis of solutions to \(
\mathbf{L}g =0 \).

\subsection*{A particular solution}
\label{sec:particular-solution}

\begin{lemma}
  \label{corollary:A.5.2}
  The particular solution to \( \mathbf{L}g=h \) with vanishing initial
  conditions (through order \( m-2 \)) is
  \begin{equation*}
    g_p(z) = \sum_{j=1}^{m-1}
    \frac{(1-z)^{-\lambda_j}}{\psi'(\lambda_j)} \int_{0}^z
    h(\zeta) (1-\zeta)^{\lambda_j+m-2} \,d\zeta.
  \end{equation*}
\end{lemma}

\subsection*{The initial conditions}
\label{sec:initial-conditions}

Having computed a basis of solutions to the homogeneous equation
and a particular solution to the inhomogeneous equation, so far we have
established equation~\eqref{1.2e} in Theorem~\ref{T:ett} modulo
determination of the coefficients \( c_1,\ldots,c_{m-1} \)
at~\eqref{eq:2.9}.

The fact that the initial conditions for \(g_p \) vanish make it
simple to solve \( \mathbf{L}g = h \) for specified initial
conditions:\ One need only match up the initial conditions of the
homogeneous solutions.

\begin{proposition}
  \label{lemma:A.2.1}
  If the general complementary solution to $\mathbf{L}g = h$ is
  written in the form
  \begin{equation}
  \label{eq:A.12.25}
  g_{\text{c}}(z) = \sum_{j=1}^{m-1} A_j (1-z)^{-\lambda_j},
\end{equation}
  then the constants \( A_j \) are given by
  \begin{equation*}
    A_j = \frac{m!}{\psi'(\lambda_j)} \sum_{k=0}^{m-2} 
    \frac{g_c^{(k)}(0)}{\rising{\lambda_j}{k+1}}, \quad j=1,\ldots,m-1.
  \end{equation*}
\end{proposition}



\section{Properties of the indicial polynomial}
\label{sec:prop-indic-polyn}

The indicial polynomial
\begin{equation}
  \label{eq:B.0}
  \psi(\lambda) \equiv \psi_m(\lambda) := \rising{\lambda}{m-1} - m!
\end{equation}
plays an important role in the analysis of random \( m \)-ary search
trees.  We will enumerate a few useful
identities involving the polynomial in this appendix.

It is well known~\cite[Chapter~3]{MR93f:68045} that \( \psi_m \) has
\(m-1\) distinct roots \( 2 = \lambda_1, \lambda_2, \ldots,
\lambda_{m-1} \) listed in nonincreasing order of real part.
As in~\cite[Chapter~3]{MR93f:68045} we introduce
\begin{equation*}
  \alpha \equiv \alpha_m:= \max_{2 \leq j \leq m-1} \Re{(\lambda_j)};
\end{equation*}
that is, \(\alpha\) is the second largest
real part among all the roots of the indicial
polynomial~\eqref{eq:B.0}.  We list some important properties of the
roots of~\eqref{eq:B.0} stated in~\cite[\S3.3]{MR93f:68045}:
\begin{enumerate}
\item The number \(-m\) is a root if and only if \(m\) is odd.  
  All other roots of \(\psi(\lambda)\) are simple, non-real roots.
\item No two roots have the same real part unless they are
  mutually conjugate.  (This follows from the strict increasingness of
  $|\rising{(s+it)}{m-1}|$ in $|t|$.)
\end{enumerate}

\subsection{Monotonicity of $\alpha$ in $m$}
\label{sec:monotonicity-alpha-m}

One can check easily that for $3 \leq m \leq 6$, $\alpha$ is strictly
increasing in $m$.  We will now prove this fact
(Theorem~\ref{theorem:alpha-monotone}) for all $m \geq 3$.  To do so we
build upon ideas in Appendix~A of~\cite{MR96j:68042}.
\begin{claim}
  \label{claim:1}
  For any $m \geq 6$ and $-\infty < x < 2$,
  \begin{equation}
    \label{eq:4}
    g_m(x) := \inf\{ y > 0: (2 - x + iy)(3 - x + iy) \dots (m - x +
    iy) \text{ is positive real} \}
  \end{equation}
  is positive and finite, and the infimum is achieved.
\end{claim}
\begin{proof}
  If $z = 2 - x + iy$ with $-\infty < x < 2$ and $y \geq 0$, then
  \begin{equation*}
    \sum_{j=2}^{m} \arg(j - x + iy) = \sum_{j=2}^m \arctan
    \frac{y}{j-x}.
  \end{equation*}
  This last expression is strictly increasing in $y \geq 0$, with
  value~0 at~0 and (for $m \geq 6$) limit 
  \begin{equation*}
    \sum_{j=2}^m \frac\pi2 = (m-1)\frac\pi2 \geq \frac{5\pi}2
  \end{equation*}
  as $y \to \infty$.  It is therefore clear that $g_m(x)$ is positive
  and finite
  and in fact is characterized as the unique root to
  \begin{equation*}
    \sum_{j=2}^m \arctan \frac{g_m(x)}{j-x} = 2\pi.
  \end{equation*}
  The set on the right in~\eqref{eq:4} is discrete, and $g_m(x)$ is
  its smallest element.
\end{proof}
\begin{claim}
  \label{claim:2}
  The function~$g_m$ in~\eqref{eq:4} is (for $m \geq 6$)
  strictly decreasing in $-\infty < x < 2$.
\end{claim}
\begin{proof}
  If $-\infty < x_1 < x_2 < 2$, then
  \begin{equation*}
    2\pi = \sum_{j=2}^m \arctan \frac{g_m(x_1)}{j-x_1} <
    \sum_{j=2}^{m} \arctan \frac{g_m(x_1)}{j-x_2},
  \end{equation*}
  so $g_m(x_2) < g_m(x_1)$.
\end{proof}
\begin{claim}
  \label{claim:3}
  The function~$g_m$ in~\eqref{eq:4} is (for $m \geq 6$)
  continuous in $-\infty < x < 2$.
\end{claim}
\begin{proof}
  Let $-\infty < x_1 < x_2 < 2$.  In Claim~\ref{claim:2} we have seen
  $g_m(x_2) < g_m(x_1)$.  To complete the proof of
  Claim~\ref{claim:3}, we will show that
  \begin{equation}
    \label{eq:10}
    g_m(x_2) \geq \frac{2-x_2}{2-x_1} g_m(x_1).
  \end{equation}
  Indeed,
  \begin{equation*}
    2\pi = \sum_{j=2}^m \arctan \frac{g_m(x_1)}{j-x_1} = \sum_{j=2}^m
    \arctan \left[\frac{j-x_2}{j-x_1} \times \frac{g_m(x_1)}{j-x_2}
    \right] \geq
    \sum_{j=2}^m \arctan \frac{\frac{2-x_2}{2-x_1}g_m(x_1)}{j-x_2},
  \end{equation*}
  so~\eqref{eq:10} follows.
\end{proof}
\begin{corollary}
  \label{corollary:4}
  For any $m \geq 6$ and $-\infty < x < 2$, define
  \begin{equation*}
    f_m(x) := (2 - x + i g_m(x))(3 - x + i g_m(x)) \dots (m - x + i g_m(x)),
  \end{equation*}
  which by Claim~\ref{claim:1} is positive real-valued.  Then $f_m$ is
  continuous.
\end{corollary}
\begin{proof}
  This is immediate from Claim~\ref{claim:3}.
\end{proof}
\begin{claim}
  \label{claim:4B}
  For fixed $-\infty < x < 2$, $g_m(x)$ is strictly decreasing in $m
  \geq 6$.
\end{claim}
\begin{proof}
  If $m \geq 6$, then
  \begin{equation*}
    2\pi = \sum_{j=2}^m \arctan \frac{g_m(x)}{j-x} < \sum_{j=2}^{m+1}
    \arctan \frac{g_m(x)}{j - x},
  \end{equation*}
  so $g_{m+1}(x) < g_m(x)$.
\end{proof}
\begin{lemma}
  \label{claim:5}
  For $m \geq 4$, let $\lambda_m = \alpha_m + i\beta_m$ (with $\beta_m
  > 0$) be a root of the indicial polynomial~$\psi_m$ with second largest real
  part.  Then $|\lambda_m + m - 1| < m + 1$.
\end{lemma}
\begin{proof}
  If $|\lambda_m + m - 1| \geq m + 1$, then by the triangle inequality
    $|\lambda_m + j| > j + 2$ for all $j=0,\dots,m-2$.
  But then
    $|\rising{\lambda_m}{m-1}| = \prod_{j=0}^{m-2} |\lambda_m + j| > m!$,
  which is a contradiction.
\end{proof}
Combining our preliminary  results we can now prove the asserted
monotonicity of~$\alpha$.
\begin{theorem}
  \label{theorem:alpha-monotone}
  The second largest real part,~$\alpha_m$, among roots of the indicial
  polynomial~$\psi_m$ is strictly increasing in $m \geq 3$.
\end{theorem}
\begin{proof}
  For $m \leq 6$ the result is easily verified.  We now show that
  $\alpha_{m+1} > \alpha_m$ for $m \geq 6$.  Observe $f_{m+1}(0) =
  |f_{m+1}(0)| > |(m+1)!| = (m+1)!$ and
  \begin{align}
    f_{m+1}(2 - \alpha_m)  &= |f_{m+1}(2-\alpha_m)|  \notag\\
    &= \bigl|[\alpha_m + i g_{m+1}(2-\alpha_m)] \dots [\alpha_m +
    (m-2) + i
    g_{m+1}(2 - \alpha_m)] \notag\\
    &\qquad {}\times[\alpha + (m-1) + i g_{m+1}(2-\alpha_m)]\bigr|
    \notag \\
    \label{eq:12}
    &< \bigl| [\alpha_m + i g_m(2-\alpha_m)] \dots [\alpha_m + (m-2) +
    i g_m(2 - \alpha_m)]\\
    &\qquad {}\times[\alpha_m + (m-1) + i g_m(2-\alpha_m)] \bigr|   \notag \\
    \label{eq:13}
    & \leq m! |\alpha_m + (m-1) + i \beta_m|\\
    &= m! |\lambda_m + m - 1| < (m+1)!.  \notag
  \end{align}
  Inequality~\eqref{eq:12} follows from Claim~\ref{claim:4B},
  inequality~\eqref{eq:13} holds since $g_m(2-\alpha_m) \leq \beta_m$,
  and the last inequality is a consequence of Lemma~\ref{claim:5}.
  Therefore, by Corollary~\ref{corollary:4}, $f_{m+1}(x) = (m+1)!$ for
  some $0 < x < 2 - \alpha_m$, and so $\rising{\lambda}{m} = (m+1)!$
  for some $\lambda$ with $\Re\lambda > \alpha_m$. That is, $\alpha_{m+1}
  > \alpha_m$, as desired.
\end{proof}
\begin{remark}[Asymptotics of $\alpha_m$]
  \label{rem:alpha-asymptotics}
  Using~(B.7) and~(B.8) from~\cite{MR90a:68012} and the
  characterization of $g_m(x)$ in Claim~\ref{claim:1}, we can
  establish
  \begin{gather*}
    \alpha_m = 2 - \left( 1 + o(1) \right)2\pi^2 \left(\frac{\pi^2}6 -
    1 \right) \ln^{-3}m, \\
    \beta_m = \left( 1 + o(1) \right) {2\pi}\ln^{-1}m;
  \end{gather*}
  we omit the proof.  Thus
  the Proposition in Appendix~B of~\cite{MR90a:68012}
  is asymptotically optimal, to first order.
\end{remark}

\subsection[Identities involving the indicial polynomial]{Identities
  involving the indicial polynomial}
\label{appendix:indicial}
\begin{identity}
  \label{identity:B.1}
  When \( \lambda \notin \{ \lambda_1,\ldots,\lambda_{m-1} \} \),
  \begin{equation*}
    \sum_{j=1}^{m-1} \frac1{(\lambda-\lambda_j)\psi'(\lambda_j)} =
    \frac{1}{\psi(\lambda)}
  \end{equation*}
\end{identity}
For \( r \) and \( n \) positive integers, let \( H_n^{(r)} \) denote
the $r$th-order harmonic number
\begin{equation*}
  \label{eq:B.4}
  H_n^{(r)} := \sum_{j=1}^n \frac{1}{j^r}.
\end{equation*}
When \( r=1 \) we will use \( H_n := H_n^{(1)} \) for the usual
(1st-order) harmonic number.

\begin{identity}
  \label{identity:B.2.2}
  For \( 0 \leq k \leq m-3 \),
  \begin{equation*}
    \sum_{j=1}^{m-1} \frac{\lambda_j^k}{\psi'(\lambda_j)} = 0.
  \end{equation*}
\end{identity}
\begin{identity}
  \label{identity:psp2}
  \[
  \psi'(2) = m! (H_m - 1) \qquad \text{and} \qquad \psi''(2) = m! [
  (H_m -   1)^2 - (H_m^{(2)} - 1)].
  \]
\end{identity}
\begin{identity}
  \label{identity:B.5.22}
  \begin{equation*}
    \sum_{j=2}^{m-1} \frac1{(\lambda_j-2)\psi'(\lambda_j)} =
    \frac1{2(m!)} \left[ 1 - \frac{H_m^{(2)}-1}{(H_m-1)^2} \right].    
  \end{equation*}
\end{identity}

\medskip
\noindent\textbf{Acknowledgment.}  The authors thank two anonymous
referees for helpful suggestions.


\bibliographystyle{habbrv}
\bibliography{msn,leftovers}

\end{document}